\documentclass[a4paper,12pt]{amsart}\usepackage{amsfonts,amssymb,oldgerm,enumerate,mathrsfs,bbm,rotating,curves,xcolor} %
\UseRawInputEncoding
\usepackage[bookmarksdepth=3]{hyperref}

\def\C{\mathbb{C}}
\def\k{\mathbbm{k}}
\def\N{\mathbb{N}}
\def\bS{\mathbb{S}}\def\R{\mathbb{R}}\def\Z{\mathbb{Z}}

\def\di{\partial}

\def\infl{\inf\limits}

\newcommand{\quots}[2]{{\footnotesize\left.\raisebox{0.4ex}{$#1$}\! / \!\raisebox{-0.4ex}{$#2$}\right.}}

\renewcommand{\stackrel}[2]{\ \lower 0.4ex \hbox{$\mathrel{\mathop{#2}\limits^{#1}}$}\ }

\def\tf{\tilde{f}}
\def\tga{{\tilde{\ga}}} 

\def\cU{\mathcal U}
\def\tf{\tilde{f}}\def\tcU{\tilde{\cU}}
\def\tX{{\tilde{X}}}\def\tY{{\tilde{Y}}}
\def\tZ{{\tilde{Z}}}

\def\hx{\hat{x}}\def\hy{\hat{y}}\def\hz{\hat{z}}

\def\al{\alpha}\def\ga{\gamma}\def\Ga{\Gamma}\def\be{\beta}\def\de{\delta}\def\De{\Delta}
\def\ep{\epsilon}\def\om{\omega}
\def\si{\sigma}\def\Si{\Sigma}

\def\cX{\mathcal X}
\def\liml{\lim\limits}

\def\uf{{\underline{f}}}\def\ug{{\underline{g}}}\def\uh{{\underline{h}}}

 \def\uom{{\underline{\om}}}\def\up{{\underline{p}}}\def\us{{\underline{s}}}

\def\ux{{\underline{x}}}\def\uz{{\underline{z}}}

\def\empty{\varnothing}

\newcommand{\bbm}{\begin{bmatrix}}\newcommand{\ebm}{\end{bmatrix}}
\newcommand{\ber}{\begin{array}{l}}\newcommand{\eer}{\end{array}}
\newcommand{\bpm}{\begin{pmatrix}}\newcommand{\epm}{\end{pmatrix}}
\newcommand{\bM}{\begin{matrix}}\newcommand{\eM}{\end{matrix}}
\newcommand{\bee}{\begin{enumerate}}\newcommand{\eee}{\end{enumerate}}
\newcommand{\bei}{\begin{itemize}}\newcommand{\eei}{\end{itemize}}

\DeclareMathOperator{\codim}{codim}
\DeclareMathOperator{\Cone}{Cone}\DeclareMathOperator{\diam}{diam}\DeclareMathOperator{\dist}{dist}
\DeclareMathOperator{\Link}{Link}\DeclareMathOperator{\Mat}{Mat}\DeclareMathOperator{\ord}{ord}

\DeclareMathOperator{\Sing}{Sing}\DeclareMathOperator{\Span}{Span}\DeclareMathOperator{\tord}{tord}
\DeclareMathOperator{\vanrate}{van.\!rate}

\def\wrt{with respect to }\def\iff{if and only if }\def\sset{\!\subset\!}\def\sseteq{\!\subseteq\!}\def\smin{\!\setminus\!}\def\sin{\!\in\!}

\def\<{\!<\!}\def\>{\!>\!}\def\={\!=\!}

\newcommand{\beq}{\begin{equation}}\newcommand{\eeq}{\end{equation}}

\newtheorem{Lemma}{Lemma}[section]\newcommand{\bel}{\begin{Lemma}}\newcommand{\eel}{\end{Lemma}}
\newtheorem{Example}[Lemma]{Example}\newcommand{\bex}{\begin{Example}\rm}\newcommand{\eex}{\end{Example}}
\newtheorem{Proposition}[Lemma]{Proposition}\newcommand{\bprop}{\begin{Proposition}}\newcommand{\eprop}{\end{Proposition}}
\newtheorem{Property}[Lemma]{Property}\newcommand{\bproperty}{\begin{Property}}\newcommand{\eproperty}{\end{Property}}
\newtheorem{Definition-Proposition}[Lemma]{Definition-Proposition}

\def\bpr{~\\{\em Proof.\ }}
\newcommand{\epr}{{\hfill\ensuremath\blacksquare}}
\newtheorem{Theorem}[Lemma]{Theorem}\newcommand{\bthe}{\begin{Theorem}}\newcommand{\ethe}{\end{Theorem}}
\newtheorem{Definition}[Lemma]{Definition}\newcommand{\bed}{\begin{Definition}}\newcommand{\eed}{\end{Definition}}
\newtheorem{Remark}[Lemma]{Remark}\newcommand{\beR}{\begin{Remark}\rm}\newcommand{\eeR}{\end{Remark}}
\newtheorem{Corollary}[Lemma]{Corollary}\newcommand{\bcor}{\begin{Corollary}}\newcommand{\ecor}{\end{Corollary}}
\newcommand{\bet}{\begin{tabular}{cccccccc}}\newcommand{\eet}{\end{tabular}}

\def\vv{{\vec{v}}}

\newcommand{\isom}{\xrightarrow[\,\smash{\raisebox{1.15ex}{\ensuremath{\scriptstyle\sim}}}\,]{}}

\title[]{F\MakeLowercase{ast vanishing cycles on perturbations of complex weighted-homogeneous complete intersection germs}}
\author[]{D\MakeLowercase{mitry}  K\MakeLowercase{erner and} R\MakeLowercase{odrigo} M\MakeLowercase{endes}}

\thanks{We were supported by the Israel Science Foundation,  grants No.  1910/18 and 1405/22}

\address{\tiny Dmitry Kerner: Department of Mathematics, Ben-Gurion University of the Negev, P.O.B. 653, Be'er Sheva 84105, Israel.
 \qquad dmitry.kerner@gmail.com}
\address{\tiny Rodrigo Mendes: Instituto de ci\^encias exatas e da natureza, Universidade de Integra\c{c}\~ao Internacional da Lusofonia Afro-Brasileira (unilab), Campus dos Palmares, Cep. 62785-000. Acarape-Ce, Brasil and  Departament of Mathematics, Ben-Gurion University of the Negev, P.O.B. 653, Be'er Sheva 84105, Israel. \quad  rodrigomendes@unilab.edu.br}
\subjclass[2020]{Primary
14J17; %Singularities of surfaces or higher-dimensional varieties
\quad Secondary
51F30 %Lipschitz and coarse geometry of metric spaces
14B05 %Singularities in algebraic geometry
32S05%Local complex singularities
}

\keywords{Singularity Theory, Lipschitz Geometry of Singularities, fast cycles, inner metrically conical germs, foliations by arcs}
\date{\today\ \  filename: \jobname.tex}

\begin{document}\setcounter{secnumdepth}{6} \setcounter{tocdepth}{1}

 \begin{abstract}
 Take a complex-analytic germ $X\sset (\C^N,o)$ with arbitrary singularity. In many cases $\Link[X] $ contains cycles that  vanish faster than linearly, when $\Link[X]$ shrinks to the origin. These ``fast cycles" capture the crucial metric/Lipschitz properties of $X.$

 We consider germs (of arbitrary dimensions and codimensions) that are perturbations of weighted-homogeneous complete intersections.
  For such germs we determine the fast cycles (i.e. their homotopy type, tangent cone, vanishing rates) via the weights. This gives a vast zoo of fast cycles with prescribed properties.

%\medskip

\noindent As an immediate application we get countable families of (distinct) exotic Lipschitz structures on germs of topological manifolds, for each fixed dimension, codimension, and multiplicity.

Another application is the obstruction for germs to be inner metrically conical, e.g.:
 \bei
 \item
 (with certain assumptions) If $X$ is IMC then the $n$ lowest weights coincide.

\item  Let the surface germ $X=V(f)\sset (\C^3,o)$ be Newton-non-degenerate and IMC.
  Then for some of the faces of the Newton diagram the two lowest weights coincide.
 \eei
 \end{abstract}
 \maketitle\vspace{-1cm}

\tableofcontents

\vspace{-1cm}

\section{Introduction}
\subsection{}\label{Sec.Intro.1.} Take a complex-analytic germ
 $X\sset (\C^N\!,o).$
 The first step in its visualization is the Conic Structure Theorem:
 {\em The germ $X$ is (subanalytically) homeomorphic to the standard cone over the link,}
  $\Cone[\Link[X]]\sset (\C^N\!,o).$
   In ``most cases" this homeomorphism cannot be chosen differentiable in any sense.

  In some cases this  topological equivalence can be strengthened to (non-embedded) bi-Lipschitz equivalence, as follows.
   The (standard) metric on $(\C^N,o)$ induces the inner metrics on $X$ and on $\Cone[\Link[X]],$ by the length of the shortest path between two points.
    The germ $X$ is called {\em inner-metrically-conical} (IMC) if the homeomorphism $X\cong \Cone[\Link[X]]$ can be chosen
    bi-Lipschitz \wrt
     these inner metrics. These IMC-germs are the simplest models from the Lipschitz point of view.

      See \cite{Handbook.IV} for a brief introduction to the Lipschitz geometry of Singularities.

\medskip

Any complex-analytic curve-germ is IMC.  Hence the inner Lipschitz geometry of these germs is trivial.

 In higher dimensions the situation is more complicated. Only sporadic examples are understood.
  The first
  weighted-homogeneous non-IMC surface germs were found in \cite{Birbrair-Fernandes.08}.
   If $X$ is a weighted-homogeneous IMC-germ of quotient type, $\quots{\C^2}{G},$ then its two lowest weights are necessarily equal, \cite{Birbrair-Fernandes.Neumann.08}.
   The converse statement was verified in \cite{Birbrair-Fernandes.Neumann.09} for  Brieskorn-Pham singularities in $(\C^3,o)$.
 Finally, \cite{Birbrair.Neumann.Pichon.14} established the thick-thin decomposition for normal surface germs.
 In particular, this gave an algorithm to verify the (non)IMC-property via the resolution graph of the singularity.
  Using this \cite{Okuma.17} proved: a Brieskorn-Pham surface germ that is an isolated complete intersection is IMC \iff its two lowest weights are equal.
 In higher dimensions, among the $A_k$-types (with equation $\sum^n_{i=1} x^2_i+x^{k+1}_{n+1}=0$), the only IMC-case is $A_1,$ \cite{Birbrair.Fernandes.Neumann.Grandjean.O'Shea14}.

\medskip

The main obstructions  to the IMC property (known up to now) are:
  fast loops, choking horns, and separating sets.

Using separating sets, W.~Neumann constructed an infinite collection of germs of topological manifolds, of any given complex dimension $\geq 3$, with pairwise distinct inner Lipschitz structures, see the appendix to\cite{Birbrair.Fernandes.Neumann.Grandjean.O'Shea14}.
 Recall that an outer bi-Lipschitz isomorphism $X\cong(\C^n,o)$ implies an analytic isomorphism, \cite{Birbrair.Fernandes.Le.Sampaio}.

 Fast loops,  and more generally, fast cycles  are important far beyond the IMC-context. They carry   essential information on the local Geometry and Topology of $X.$

\medskip

To detect or to rule out  fast loops is in general a hard task. Besides particular examples (with complicated proofs),
 the only `manageable' cases are the normal surface singularities. There the  criteria via   resolution graphs are known,   \cite{Birbrair.Neumann.Pichon.14}.
  In \cite{Kerner-Mendes.Coverings.Discriminants} we give another approach, detecting fast loops via ramified coverings and their discriminants.

\medskip

 Now we  show how  to detect fast cycles on  perturbations of weighted-homogeneous germs.
 Moreover, we describe their homotopy type, vanishing rates, and tangent cones.

\medskip

Our criteria for fast cycles can probably be stated in terms of the general theories of \cite{Birb.-Brass.00}, \cite{Valette}, \cite{Bob.-Hein.-Per.}, \cite{Bob.-Hein.-Per.-Sam.}. However, it is not clear whether these theories allow one to get such criteria, or at least to simplify the proofs.

\subsection{Generalities on fast cycles}
Fast loops were defined and studied in particular cases, see \S\ref{Sec.Intro.1.}. But we do not know a general reference on fast cycles, especially when their  tangent cones are of positive dimensions. Therefore we introduce the general notions in \S\ref{Sec.Fast.Cycles}, in particular define the homotopy type of the fast cycle, and the vanishing rate. We give many examples and prove:

{\bf Lemma  \ref{Thm.IMC.germs.have.no.fast.cycles}.}
 {\em A (subanalytic, closed) inner-metrically-conical germ  has no fast cycles.}

Therefore fast cycles are obstructions to the IMC property.

\subsection{Fast cycles on germs with (perturbed) weighted-homogeneous foliation}   \

\noindent\parbox{13.5cm}{Suppose a subanalytic germ $X$ is foliated by (real) arcs that are ``weighted-homogeneous modulo higher order terms".
  (The simplest example is just a weighted-homogeneous germ.)
  Then any cycle on the link, $Z \sset  \Link[X],$   defines the subgerm $\R_{\ge0}Z \sset  X,$ the union of the arcs through $Z.$

  We give a very general Lemma \ref{Thm.IMC.restriction.weights.general}, stating the minimal assumptions to ensure  that $ Z$ is a fast cycle.
 Moreover, the tangent cone and the vanishing rate of $\R_{\ge0}Z$ are controlled by the weights.}
\begin{picture}(0,0)(-17,30)
\color{blue}\thicklines
%\curve(0,0,10,0,20,2,30,6,40,12,50,20,60,35)\curve(0,0,5,1,10,3,15,6,20,10,25,15,30,23,35,35,40,49)
\curve(0,0,10,0,20,2,30,6,40,12,50,20)\curve(0,0,5,1,10,3,15,6,20,10,25,15,30,23,35,35)

\curve(52,23,54,26) \curve(56,29,58,32)
\curve(36,38,37,41) \curve(38,44,39,47)

\curve(44,15,33,15.5,28,19)\curve(38,11,29,11.5,25,15)
\curve(32,7,25,7.5,21,11)\curve(22,3,17,3.5,16,7)

\curve(16,1,13,1.5,12,4)

\put(0,18){$\scriptstyle{\pmb{\R_{\ge0}Z}}$}

\color{red}
\curve(60,35,45,37,40,49)\curve(60,35,55,45,40,49)  \put(43,53){$\scriptstyle{\pmb Z}$}

\curve(43,46,45,42,50,39)\curve(42,45,47,44,49,38.5)

\curve(51,43,53,39,58,36)\curve(50,42,55,41,57,35.5)
\thinlines
\color{black}
\curve(20,40,40,30,50,10)\curve(20,40,35,50,45,70)\curve(45,70,60,50,80,40)\curve(50,10,60,30,80,40)

% \curve(45,65,60,48,75,40) \curve(45,60,55,48,70,40)

\put(50,5){$\scriptstyle{\pmb{ \Link[X]}}$}

\end{picture}

\subsection{Fast cycles on perturbations of weighted-homogeneous germs}\label{Sec.Intro.3.}

 Take a complex  weighted-homogeneous complete intersection germ $X_o:=V(\uf_\up)\sset (\C^N,o)$
  of dimension $n.$  (See \S\ref{Sec.Notations.Basics} for the notation.) It can have arbitrary singularity, e.g. can have multiple components.
  Suppose not all of the coordinate-weights are equal, though some of them might coincide. Thus we arrange:
   \beq\label{Eq.weights.split}
   \om_1=\cdots=\om_r< \om_{r+1}\le\cdots\le\om_N, \quad
   \text{for some } 1\le r<N.
   \eeq
Perturb $\uf_\up$ by higher-order terms, $X\!:=V(\uf_\up\!+ \uf_{>\up})\sset (\C^N,o),$ here $ \ord_\om f_{p_i} \!<\! \ord_\om f_{>p_i} .$

In the   surface case, $n=2,$    we have (under certain assumptions):
\\{\bf Proposition \ref{Thm.Fast.Cycles.via.weights.Surfaces}.} {\em   If  $\om_1<\om_2$ then $X$ has a fast loop of vanishing rate$\ge\frac{\om_2}{\om_1}.$}

 \medskip

\noindent In arbitrary dimension  denote $l\!:=\!n-\dim_\C \Sing[X\!\cap\! V(x_1)].$ Suppose $\om_1\!<\!\om_l,$  thus $1\!\le\! r \!<\!l\!\le\! n.$
\\{\bf Theorem \ref{Thm.IMC.criterion.via.weights}.} (under certain assumptions) {\em
        Then $X$ has a fast cycle $\R_{\ge0}Z,$  of   dimension $l \! \le\! \dim_\R \R_{\ge0}Z \!\le\! n ,$ and with the properties:
\bee
\item    Its tangent cone   lies in the half-plane $\R_{\ge0}\hx_1\!\times\! \C^{r-1}_{x_2\dots x_r},$
   and  is of dimension $\dim_\R T_{\R_{\ge0}Z}\!\le\! r.$
  \item     The hornic homotopy type of $Z$ is that of the Milnor fibre of the germ $X\cap V(x_1).$
\item  The vanishing rate is bounded, $\vanrate_X(\R_{\ge0}Z)\ge\frac{\om_{r+1}}{\om_1}.$

\item  Under additional assumptions the vanishing rate equals $\frac{\om_{r+1}}{\om_1}.$
  \eee
}   %Therefore, if $X$ is IMC, then $\om_1=\cdots=\om_l$.

\medskip

With this result the detection of fast cycles becomes immediate, and we get the whole zoo of them, with various invariants.
 E.g. for (semi)weighted-homogeneous ICIS we have.
\\
{\bf Corollary \ref{Thm.IMC.criterion.semi.weighted.homogen.germs}.}
 {\em Suppose the germs $X_o$ and  $X\!\cap \!V(x_1)$ are ICIS. Take the Milnor number $\mu\!:=\!\mu(X\!\cap\! V(x_1)).$
\bee

\item If $\om_1\!<\!\om_n$ then $X$ has a fast cycle of homotopy type
  $\vee^{\mu}S^{n-1} $
 and of the vanishing rate at least  $ \frac{\om_{r+1}}{\om_1} .$  With a simple additional assumption the vanishing rate is exactly  $ \frac{\om_{r+1}}{\om_1} .$

The tangent cone of this fast cycle lies in the half-plane $\R_{\ge0}\hx_1\times \Span_\C(\hx_2,\dots,\hx_{r})$ and is of $\dim_\R\!\le\!r.$
\item
  If $X$ is IMC  then $\om_1=\cdots=\om_n.$
\eee}
\noindent Immediate examples are Brieskorn-Pham hypersurface germs, and Brieskorn-Pham complete intersections (in all dimensions and codimensions), and their perturbations, see Example \ref{Ex.Examples.of.Fast.Cycles}.
\\
Theorem \ref{Thm.IMC.criterion.via.weights} applies also when the germ $X_o$ has a non-isolated singularity, see Corollary \ref{Thm.fast.cycles.on.perturb.weighted.homogen.hypersurf.sing}.

\subsection{Applications}
We get a simple way to cook various exotic Lipschitz structures. E.g. in \S\ref{Ex.Exotic.Spheres} we give multi-index collections of Lipschitz-distinct topological manifolds.
\\{\bf Corollary \ref{Thm.Exotic.Structures}.} {\em
     For each triple $(n,p,r),$ with $n\!\ge\!3,$ $p\!\ge\!2,$   $r\!<\!n,$ there exists a countable number of weighted-homogeneous hypersurface-germs $X\sset (\C^{n+1},o)$, each having $mult(X)\!=\!p,$ and weights $\om_1\!=\!\cdots\!=\!\om_r\!<\!\cdots,$
 all homeomorphic to $(\R^{2n},o),$ but   of (pairwise) distinct inner-Lipschitz-types.
    }

Recall, that in \cite{Birbrair.Fernandes.Neumann.Grandjean.O'Shea14} the Lipschitz structures were distinguished by (arbitrary large number of) separating sets. We distinguish these structures by (arbitrary large) vanishing rates of fast cycles or by their topological complexity.

It is remarkable that these exotic structures are realized as (very simple) complex-algebraic hypersurfaces.

\medskip

\noindent As another application, take a Newton-non-degenerate power series $f\sin \C\{x_1,\dots,x_{n+1}\}.$ Fix a face of the Newton diagram, $\si\sset \Ga_f\sset \R^{n+1}.$ Consider $f$ as a perturbation of its restriction to $\si,$ i.e. $f=f_{|\si}+{f_{>\si}}.$ This allows to detect fast vanishing cycles just from the weights of $\si.$ Vice-versa, suppose $V(f)\sset (\C^{n+1},o)$ is IMC. Then the possible weights of $\si$ are restricted.
\\{\bf Example \ref{Ex.IMC.that.are.NND}.} Suppose the face  $\si$   intersects a coordinate hyperplane of the lowest weight.
\bei
\item ($n=2$, surfaces) Then the two lowest weights of $\si$ coincide.
\item ($n=3$, three-folds) Suppose the distance of $\si$ from each coordinate hyperplane (in $\R^4$) is at most 1. Then the three lowest weights of $\si$ coincide.
\eei

\subsection{} Our construction of fast cycles (e.g. the proof of Theorem \ref{Thm.IMC.criterion.via.weights}) gives much stronger invariants, e.g. the monodromy and intersection lattices of fast cycles. But this is postponed to the next paper.

\subsection{Declarations}
The author declare that they have no conflict of interest. There are no data sets generated during the current
study.

\subsection{Acknowledgements} Our thanks are to Lev Birbrair, Alexandre Fernandes, Andrey Gabrielov, Edson Sampaio, for important advices, and to the referee for the careful reading.

 Most of the results were obtained (and announced) back in 2022-2023, but the finalizing step was delayed by the war. D.K. thanks BGU for continuous support.

\section{Notations, conventions, preliminaries}\label{Sec.Notations}
\noindent Let $\k=\R,\C. $ We recall some  basics of Lipschitz geometry of singularities, see \cite{Handbook.IV}.

\makeatletter
\renewcommand{\p@enumi}{\thesection.}
\makeatother
\subsection{}\label{Sec.Notations.Basics}
\bee[\bf i.]

\item\label{Sec.Notations.Basics.i}
Through the paper by germs we mean their small representatives. The base point of a germ is always $o\in \k^N$ and is usually omitted.
All our germs are subanalytic and (unless stated otherwise) closed. %And their links as well.

A  neighborhood of a subgerm $Y\sset X$ is a (subanalytic) germ $Y\sseteq  \cU(Y)\sseteq \! X$ such that $\cU(Y)\smin o$ is open inside $X.$
\\
By $\hx,\hy,\hz,$ or $\hx_i$ we denote the unit vectors along the corresponding coordinates axes in $\k^N.$

\item\label{Sec.Notations.Basics.ii}
The tangent cone $T_{X}\sset \R^N$ of the germ $X\sset(\R^N,o)$ is defined (set-theoretically) as the Whitney cone $C_3$,
  \cite[pg.210]{Whitney},  \cite[Chapter 2]{Chirka}.
  The tangent cone is subanalytic as well, with $\dim_\R T_{X} \le \dim_\R {X}.$ In particular, if  $\dim_\R T_{X}=1,$ then $T_{X}$ is a finite union of lines and half-lines.

While $T_X$  depends essentially on the embedding $X\sset (\R^N,o),$ not just on the Lipschitz geometry of $X,$  $\dim_\R T_X$ is determined by the inner Lipschitz type of $X.$ See \cite{Sampaio}.

\item\label{Sec.Notations.Basics.v} All our arcs are subanalytic germs diffeomorphic to $(\R_{\ge0},o).$ Usually we take the  parametrization  $\ga_\us(t)=(s_1\cdot t+o(t),\dots,s_N\cdot  t+o(t)),$ with $\us\in\bS.$
 This coincides with length-parametrization up to terms of order $o(t)$.
\\
Take a subanalytic function-germ $X\stackrel{f}{\to}\R.$ Its restriction to an arc is either identically zero or $f(\ga(t))=c\cdot t^\al+o( t^\al)$ with $c\neq 0.$ Accordingly define $\ord(f|_\ga)=\infty$ or $\ord(f|_\ga)=\al.$

A foliation of the germ $X$  by a set of arcs  is  a subanalytic family $\{\ga_\us\}_{\us\in \Link[X]}$ that covers a representative of $X,$
   and satisfies (for this representative): $\ga_\us\cap \ga_{\us'}=o$ for $\us\neq\us'.$

\item\label{Sec.Notations.Basics.iii}
Take the unit sphere   $\bS:=\{ x|\ \| x\|=1\}\sset \k^N.$  This is $ S^{N-1}\sset \R^N$ or $ S^{2N-1}\sset \C^N.$

We use multi-indices, $\ux=(x_1,\dots,x_N)\in \k^N,$  $\us=(s_1,\dots,s_N)\sset \bS,$ $\uom=(\om_1,\dots,\om_N),$ $\up=(p_1,\dots,p_c),$ $\uf_\up=(f_{p_1},\dots,f_{p_c}).$
 Accordingly  we abbreviate:
 \beq
 t^{\uom}\cdot \us\!: =\!(t^{\om_1}  s_1 ,\dots,t^{\om_N} s_N  ), \
 t^{\up}\cdot \uf_\up\!: =\!(t^{p_1}  f_{p_1} ,\dots,t^{p_c}  f_{p_c}  ), \
  t^{-\up}\cdot \uf_{>\up}\!\!: =\!(t^{-p_1}  f_{\!>p_1} ,\dots,t^{-p_c}  f_{\!>p_c}  ).
 \eeq

Here $f_{p_i}\!\neq\! 0$ is weighted-homogeneous of weight $p_i,$ i.e. $f_{p_i}(t^{\om_1}x_1,\dots,t^{\om_N}x_N)=t^{p_i}\cdot f_{p_i}(x).$
 Therefore $\ord ( f_{p_i}\big|_{t^{\uom}\cdot\us})= p_i $ and $\ord ( f_{>p_i}\big|_{t^{\uom}\cdot\us})>p_i.$
\item\label{Sec.Notations.Basics.iv}
Occasionally we pass to the polar coordinates, $\bS\times\R_{\ge0}\stackrel{\si}{\to}\k^N,$ by $(\us,t)\to t\cdot \us.$

  Any subset $X\sset \k^N$ lifts to its strict transform, $\tX:=\overline{\si^{-1}(X\smin o)}\sset\bS\times\R_{\ge0}.$ If $X$ is closed, resp. compact, resp. (sub)analytic, then  $\tX$ is of the same kind.

Every function  $f:(\k^N,o)\to(\k^1,o)$ lifts to its pullback    $\tf:=f\circ \si: \bS\times \R_{\ge0}\to (\k^1,o).$

 The link at radius zero is $\Link_o[X]:=\tX\cap (\bS\times\{o\}).$ Observe:
 \beq
 \Link_o[X]=\Link_o[T_X] \qquad \text{and}\qquad T_X=\si(\Link_o[X]\times \R_{\ge0}).
\eeq

\medskip

The embedded topological type of the link of a germ, $\Link_\epsilon[X]=X\cap \bS_\ep\sset \bS_\ep,$ is well-defined.
 Moreover, the Lipschitz type is well-defined, \cite{Valette.7}.

\item\label{Sec.Notations.Basics.vi}
An embedded subgerm $X\sset (\R^N,o)$ acquires the outer metric, $d_{(\R^N,o)}(x,x'):=\|x-x'\|.$ In addition it gets the inner metric,
  $d_X(x,x'),$ defined as the infimum of lengths of paths $[x\rightsquigarrow x']$ inside $X.$
  E.g. for any convex subset $X\sseteq\R^N$ this is the usual Euclidean distance,   $d_X(x,x')=\|x-x'\|.$

   This function, $X\times X\stackrel{d_X}{\to}\R_{\ge0}$ can be pathological in various ways.
    However, by \cite{Kurdyka.Orro.97} there exists another,   subanalytic, metric  $d^{sub}_X,$ that is equivalent to $d_X.$
 Throughout the paper, we always use $d^{sub}_X.$

 A subgerm   $X\sset (\R^N,o)$ is called Lipschitz normally embedded (LNE), if the outer and inner metrics are equivalent, i.e. $0<c\le \frac{d_X(x,x')}{d_{(\R^N,o)}(x,x')}\le C<\infty$ for all $x\neq x'\in X.$

\item\label{Sec.Notations.Basics.vii}
  The (inner) tangency order of two arcs $\ga_1, \ga_2 \sset X$ is $\tord_X (\ga_1,\ga_2):=\ord_t d^{sub}_X(\ga_1(t),\ga_2(t)),$
 with length parametrization of $\ga_1,\ga_2$.
 Note that this is independent of the choice of $d^{sub}_X$ approximating $d_X.$

The diameter of a subgerm $Y\sset X$ is the lowest tangency order of its arcs,
\beq
\diam_X[Y]:=\inf_{\ga,\ga'\sset Y}  \tord_X(\ga,\ga') .
\eeq
If $X\sseteq(\R^N,o)$ and $\Link[Y]$ is connected, then the inner and outer diameters coincide,
 $\diam_X Y=\diam_{(\R^N,o)}Y.$
  This follows by the pancake decomposition, \cite{Birbrair.Mostowski.00}.

\item\label{Sec.Notations.Basics.viii} Let $X\sset \R^N$ be a sub-analytic set.
Every homology class in $H_*(X,\Z)$ admits a subanalytic representative.
  Indeed, the singular homology $H_*(X,\Z)$ coincides with the simplicial homology.
 And sub-analytic manifolds admit subanalytic triangulations, \cite[\S4.3]{Coste}.

\item\label{Sec.Notations.Basics.ix} Take a $\C$-analytic complete intersection $X=V(f_1,\dots,f_c)\sset (\C^N,o).$
 Here $c<N$ and $\dim_\C X=N-c.$ Take the  Jacobian matrix $[f']:=[f'_1,\dots,f'_c]^T\in \Mat_{c\times N}.$
  The singular locus $\Sing(X)$
 is defined (inside $X$) by the ideal of size-$c$ minors $I_c[f'].$
 Thus $\Sing(X)\sseteq X$ is a $\C$-analytic subgerm.
\\
\parbox{11cm}{

\item\label{Sec.Notations.Basics.x} (Embedded conic structure theorem, e.g. theorem 9.3.6, page 225 of \cite{BCR}) Every subanalytic germ can be radially rectified,
 as on the diagram. Here $\varphi$ is a subanalytic homeomorphism that acts linkwise, i.e. $\varphi(\bS_\ep)=\bS_\ep$ for all $0\le \ep\ll1.$}
\qquad $\bM (\R^N,o)\stackrel{\varphi}{\isom}(\R^N,o)\hspace{0.5cm}\vspace{-0.2cm}\\\cup \hspace{2cm}\cup\hspace{1cm}\\ X\isom{} \Cone[\Link[X]]\eM$
\eee

\subsection{Hausdorff arc-distance}\label{Sec.Notations.Hausdorff.Distance}
The  tangency order of  an  arc and a   subgerm,  $\ga , Y\sset X,$ is $\tord_X(\ga,Y\!)\!:= \!\sup_{\!\ga_Y\sset Y} \tord_X(\ga,\ga_Y).$
 For closed subanalytic germs this supremum is realized as the maximum.

The ``Hausdorff arc-distance" between (closed, subanalytic) subgerms,   $Y_1 , Y_2\sset X,$ is
\[
\dist_X(Y_1,Y_2):=\min\Big\{\infl_{\ga_1\sset Y_1}\{\tord_X(\ga_1,Y_2)\},\quad
  \infl_{\ga_2\sset Y_2}\{\tord_X(Y_1,\ga_2)\}\Big\}\le \infty.
 \]
Both infima are realized as minima.  This implies the basic properties of the distance:
\bee[\bf i.]
\item $\dist_X(Y,Y')=\infty$ \iff $Y=Y'.$
\item The ultrametric triangle inequality, for any germ $Z\sset X:$

\hspace{2cm}$\dist_X(Y_1,Y_2)\ge \min\big\{ \dist_X(Y_1 ,Z),\dist_X(Z,Y_2)\big\}.$

\item (an implication for tangent cones) Take the germs $Y,Y'\sset X\sseteq (\R^N,o)$ and their tangent cones, $T_Y,T_{Y'}\sset \R^N.$
  If $\dist_X(Y,Y')>1$ then $T_Y=T_{Y'}.$
\\Vice-versa, if $T_Y\=T_{Y'}$ and $X\sseteq (\R^N,o)$ is LNE (see \S\ref{Sec.Notations.Basics.vi}), then
$\dist_X(Y,Y')\>1.$

\item Given a subgerm $Y\sset X,$ suppose $T_Y\sset X.$
Take a common neighborhood of a germ and of its tangent cone, $Y,T_Y\sseteq \cU\sseteq X.$
 Then (obviously):
 \beq
 \dist_X(Y,T_Y)\ge \dist_X\big( \overline{\cU(Y)},T_Y\big) \qquad \text{ and } \qquad \dist_X(Y,T_Y)\ge \dist_X\big(Y,\overline{\cU(Y)}\big).
 \eeq

 \medskip

 Moreover,  we claim: $\dist_X(T_Y,\overline{\cU(Y)})=\dist_X(Y,\overline{\cU(Y)}).$
\bpr  We have   $\dist_X(T_Y,\overline{\cU(Y)})\ge \min\{\dist_X(T_Y,Y),\dist_X( Y,\overline{\cU(Y)})\}=\dist_X( Y,\overline{\cU(Y)})$ and
$\dist_X( Y,\overline{\cU(Y)})\ge \min\{\dist_X(Y,T_Y),\dist_X(T_Y,\overline{\cU(Y)})\}=\dist_X(T_Y,\overline{\cU(Y)}).$
 \epr
\eee

\subsection{Homologies of the Milnor fibre}\label{Sec.Notations.Milnor.Fibre.Homology}
Take a reduced complete intersection germ $X\sset (\C^N,o),$ of $\dim_\C X=n,$ possibly with non-isolated singularity.
 Denote   $l\!:=\!n-\dim_\C \Sing[X],$ thus $1\le l\le n.$ Then the Milnor fibre    $X_F$   is $(l-1)$-connected. (See \cite{Kato.Matsumoto} for hypersurface germs and  \cite{Greuel} for complete intersections.)  In particular, $X_F$ is connected.

Moreover, the   complex analytic manifold $X_F$ of $\dim_\C X_F\!=\!n$ admits deformation-retraction to a subset of $\dim_\R\le n.$
  E.g. by the theorem of Andreotti-Frankel, \cite{Hamm-Le}.

  Therefore the only possibly non-vanishing (reduced) homologies are $H_l(X_F,\Z),$ \dots, $H_n(X_F,\Z).$
  \bel\label{Thm.Milnor.Fibre.non.vanish.homol}
  If $X$ is singular, then the Milnor fibre is non-contractible. Moreover,
  \bei
  \item either $\pi_1(X_F)\neq0,$ thus   $l=1;$
  \item or $\pi_1(X_F)=0,$  but at least one of $H_l(X_F,\Z),$ \dots, $H_n(X_F,\Z)$ does not vanish.
  \eei
  \eel
\noindent For complete intersections with isolated singularities this statement is trivial.
\\\bpr
 {\em ($X_F$ is non-contractible)}  Fix some defining equations $X=V(f_1,\dots,f_c)\sset (\C^N,o),$ where $c=\codim_\C X.$
  Take the auxiliary germ $X'=V(f_2,\dots,f_c)\sseteq (\C^N,o).$
  \bei
  \item $X'$ is a complete intersection. Otherwise $X'$ contains a component of $\dim_\C>n+1$. But then $X=X'\cap V(f_1)$ contains a component of $dim>n$. Hence the contradiction with $\dim_\C X=n.$
\item $X'$ is a reduced space. Otherwise $X'$ contains a multiple component, say $X'\supseteq Y.$ But then $X$ contains a multiple component, $X\supset Y\cap V(f_1).$     Hence the contradiction with $X$ being reduced.
  \eei
If the equations are chosen generically, then $X_F:=V(f_1-t_o,f_2,\dots, f_c)$ is the Milnor fibre of $X.$ Therefore $X_F$ can be presented as the Milnor fibre for the function $X'\stackrel{f_1}{\to}(\C^1,o).$ The latter is the Milnor fibre of a reduced hypersurface-germ inside a reduced analytic space.    For such a hypersurface singularity the total trace of monodromy (on all the vanishing homologies of $X_F$) vanishes, by  \cite[Th\'eor\`{e}me 1 bis]{A'Campo}.
  But $X_F$ is connected, thus $H_0(X_F,\Z)=\Z.$ And thus the monodromy on $H_0$ is identity. Therefore, the vanishing of the
   total trace implies: at least one other homology is non-zero.
 And thus $X\cap V(x_1-t_o)$ is non-contractible.

  {\em  (Non-vanishing homology)} Suppose $\pi_1(X_F)=0 $ and  $H_l(X_F,\Z)=0=\cdots =H_n(X_F,\Z).$
  Thus  $H_j(X_F,\Z)=0$ for all $j\ge1.$
  By Hurewicz theorem we have  $ \pi_j(X_F)\isom{} H_j(X_F,\Z) $ for all $j\ge 2.$ (See e.g. Theorems 4.32 and 4.37 of \cite{Hatcher}.)

 Therefore  $ \pi_j(X_F)=0$ for all $j\ge 1.$
 But then, by Whitehead's theorem, e.g. Theorem 4.5 of \cite{Hatcher}, $X_F$ must be contractible.
  Which is a contradiction.
 \epr

\section{Fast cycles}\label{Sec.Fast.Cycles}
  A fast cycle is   a thin germ that does not admit  hornic link-homotopy onto a non-thin germ.
 Fast loops have been studied in several works, see \S\ref{Sec.Intro.1.}, but (it seems) fast cycles of higher dimensions were never properly defined.

\subsection{Basic notions} Fix a  subgerm $Y\sset X,$ subanalytic, closed, not necessarily pure-dimensional, possibly with non-connected link.
\bed\label{Def.Fast.Cycles}
 \bee
\item  $Y$  is called {\em thin} if $\dim_\R Y> \dim_\R T_{Y}.$
\item
A subanalytic neighborhood $Y\sseteq \cU(Y)\sseteq X$ is called {\rm hornic}  if $\dist_X(Y,\overline{\cU(Y)})>1.$
\item A link-homotopy of $Y\sset X$ is a continuous map $\Phi\!\!: Y\times[0,1]\to X,$ denoted  $\Phi_t,$ satisfying:
\[
\Phi_o=Id_Y, \quad \quad\quad
\Phi_t(\Link_\ep[Y] )\sseteq \Link_\ep[X] \quad  \text{ for all $t\in[0,1]$\quad  \text{ and all}\quad   $0<\ep\ll1.$}
\]
\item A link-homotopy is called hornic if $\Phi(Y\!\times\![0,1])\sseteq  \cU(Y)\sseteq X$ for a hornic neighborhood $\cU(Y).$

\item A thin subgerm $Y\sset X$
 is called a \underline{fast cycle} if  $Y$  admits no
 hornic link-homotopy to a non-thin germ.
\quad
A fast cycle whose link is $S^1$ is called \underline{a fast loop}.

\item
A fast cycle $Y_o\sset X$ is called ``of homotopy type $Z\sset \Link[X]$", if there exists a hornic link-homotopy $\{Y_t\}$ (inside $X$) with $\Link[Y_1]=Z.$
  \eee
\eed
\noindent\parbox{12.5cm}{The picture shows a thin germ $Y_1,$ for $X=\R^3,$ with topologically non-trivial $\Link[Y_1],$ and a hornic homotopy $Y_t$ to a non-thin germ $Y_0$.

The notions of this definition depend  on the (subanalytic) inner-Lipschitz type of $X$ only, not on a particular embedding $X\hookrightarrow (\R^N,o).$
 (See \S\ref{Sec.Notations.Basics.ii} for the invariance of $\dim_\R T_Y$.)
}\quad \begin{picture}(0,0)(-25,30)
 \color{blue}
 \curve(0,0,4,4,8,9,12,14,16,22,18,27,20,32,22,40,22,45) \curve(0,0,4,4,7,9,10,14,12,20,12,28,10,38,7,50)
 \curve(22,45,12,46,7,50) \curve(22,45,18,48,7,50)
 \put(0,54){$\scriptstyle{\Link[Y_1]}$}   \put(0,17){$\scriptstyle{Y_1}$}

 \color{black}
\curve(0,0,40,40)  \curve(0,0,4,4,8,9,12,14,16,20,20,27,24,35,28,44)
\curve(28,44,32,41,40,40)\curve(28,44,34,43,40,40)
 \put(36,44){$\scriptstyle{\Link[Y_t]}$}  \put(26,34){$\scriptstyle{Y_t}$}

 \color{red}
\curve(4,4,8,7,12,10,16,13,20,16,24,18,28,20,32,22,52,32)\put(50.5,31.7){$.$}\put(50.5,31.3){$.$}
\put(54,28){$\scriptstyle{\Link[Y_0]}$}   \put(14,5){$\scriptstyle{Y_0}$}
\end{picture}

\bex\label{Ex.Fast.Cycles}
%\leavevmode\vspace{0.4cm}
\bee[\bf i.]
\item Take the $\be$-horn, $Y:=\{x^2+y^2=z^{2\be},\ z\ge0\}\sset (\R^3,o),$ with $\be>1.$ Then $T_Y=\R_{\ge0}\hz$ and $Y$ is thin.
 $Y$ is hornically link-homotopic to $T_Y,$ if considered as a subgerm of $(\R^3,o).$
 But $Y$ admits no hornic link-homotopy to a non-thin germ if considered inside
 $X:=\R^3\smin \R_{\ge0}\hz.$    Hence $Y\sset X$ is a fast loop.
  Note that $\Link[Y]$ is contractible inside $\Link[X],$ though not hornically.
   Observe that $\diam(Y)=\be.$ This is called the vanishing rate in \S\ref{Sec.Preliminaries.Vanishing.Rate}.

In Example \ref{Ex.Vanishing.Rate} we extend this to fast cycles with $\dim_\R T_Y>1.$

\item
 More generally, take an arc and its hornic neighborhood, $\ga\sset \cU(\ga)\sset X.$ If $\cU(\ga)$ is hornically link-homotopic to $\ga,$ then $\cU(\ga)$ has no fast cycles.
 But a non-hornic link-homotopy $\cU(\ga)\rightsquigarrow \ga$ does not prevent fast cycles.

\item Let $Y$ be a fast cycle with $\dim_\R T_Y=1.$ Then $T_Y$ is a (finite) union of half-lines.
 In particular, the tangent cone of a fast loop is   a half-line.

\item A germ can be non-thin and yet possess fast cycles. E.g. $\{x^2+y^2\le z^2\}\smin \{x^2+y^2< z^{2\be}\}$
 or $\{x^2 \le z^2\}\smin \{x^2+y^2< z^{2\be}\}.$

\item If a (subanalytic) neighborhood  $Y\sset \overline{\cU(Y)} $ is hornic, then  $T_{Y}=T_{ \overline{\cU(Y)} }.$ In particular, if $\dim_\R\cU(Y)> \dim_\R Y,$ then $\overline{\cU(Y)}$ is a thin set.
\\Vice-versa, if a neighborhood $\overline{\cU(Y)}\sset X$ is LNE, see \S\ref{Sec.Notations.Basics.vi}  and $T_Y\!=\!T_{\overline{\cU(Y)}},$ then $\overline{\cU(Y)}$ is  hornic,
  see \S\ref{Sec.Notations.Hausdorff.Distance}.iii.

   \item Suppose $\cU$ is a common neighborhood of a germ and of its tangent cone, $Y\cup T_Y\sseteq \cU\sseteq X.$
 Then $\cU$ is a hornic neighborhood of $Y$ \iff $\cU$ is a hornic neighborhood of $T_Y.$
 (See \S\ref{Sec.Notations.Hausdorff.Distance}.iv.)

\item
Suppose $\dim_\R X\ge2,$ $ T_X$=half-line, and $\Link[X]$ is connected and non-contractible. Then $X$ is a fast cycle.
  \item
Suppose  $\Link[X]$ is topologically non-trivial. (E.g. $H_j(\Link[X],\Z)\!\neq\!0$ for some
 $0\!<\!j\!<\!\dim_\R \Link[X].$)
 Then one can expect fast cycles on $X.$
  However, fast cycles exist also on exotic spheres, where  $\Link[X]\cong S^{2n-1},$ see \S\ref{Ex.Exotic.Spheres}.
\eee
\eex

\subsection{Fast cycles are obstructions to IMC}
\bel\label{Thm.IMC.germs.have.no.fast.cycles}
A (subanalytic, closed) inner-metrically-conical germ  has no fast cycles.
\eel
\bpr  (Re)Embed the IMC germ as a cone, $X=\Cone[\Link[X]]\sset \R^N.$ Thus $T_X=X.$ W.l.o.g. we assume: $\Link[X]$ is connected.
 Take any (subanalytic) thin set $Y\sset X.$ We construct a hornic link-homotopy $Y\!\rightsquigarrow\! T_Y$ inside $X.$ Define two flows.
\bei
\item The radial flow, $\Phi^X_\tau\!: X\!\to\! X,$ with $\tau\!\in \![0,\infty),$ is defined by $\vv\to \tau\cdot \vv.$ Thus $\Phi^X_\tau$ preserves the rays through the origin, and sends links to links, $\Phi^X_\tau \Link_t [X]=\Link_{\tau\cdot t}[X].$
\item Take a flow contracting $Y$ to the origin and compatible with link-structure. Namely: a family of (subanalytic) continuous maps $\Phi^Y_\tau: Y\to Y,$ $\tau\in [0,1],$ satisfying:
\beq
\Phi^Y_1\!=\!Id_Y, \  \Phi^Y_0(Y)\!=\!o, \quad
 \Phi^Y_\tau[\Link_t Y]\!=\!\Link_{\tau\cdot t}[Y],  \quad     \Phi^Y_\tau \text{ is a homeomorphism for }\tau\!>\!0.
\eeq
This $\Phi^Y_t$ is defined, e.g. via the embedded conic structure theorem, \ref{Sec.Notations.Basics.x}, as $\varphi^{-1}$  of the radial contraction of $\varphi(Y).$
\eei

\medskip

Finally, define the homotopy  $\Psi\!:Y\!\times\!(0,1]\to X$   by $\Psi_\tau(y)\!=\!\Phi^X_{\frac{1}{\tau}}( \Phi^Y_\tau (y)).$
  This is a family of  continuous maps (for $\tau\!>\!0$). It acts linkwise, $\Link_t[Y]$ is sent into $\Link_t[X].$

   Pass to the polar coordinates (see \S\ref{Sec.Notations.Basics.iv}), to get the strict transforms $\widetilde{\Psi_\tau(Y)}\sset \bS\times\R_{\ge0}.$  This is a subanalytic family, satisfying:
    \bei
    \item $\widetilde{\Psi_\tau }|_{t=0}\=Id_\bS$ for all $\tau\>0.$ In particular,
     this homotopy preserves $T_Y,$ i.e.
    $T_{ \Psi_\tau Y}\=T_Y$ for each $\tau\>0.$

    \item   The limit   is a cone, $\liml_{\tau\to0}(\Psi_\tau Y)\!=\!T_Y.$ \quad
     Indeed,  $\liml_{\tau\to0}(\widetilde{\Psi_\tau(Y)})\=\Link_o[Y]\times \R_{\ge0}=\widetilde{T_Y}.$
    \eei

\noindent
 We have constructed the link-homotopy  $Y\!\stackrel{\Psi}{\rightsquigarrow}\! T_Y.$ It occurs inside the neighborhood $Y\sseteq  \Psi(Y\!\times\![0,1])\sseteq  X.$
 And $T_{\Psi(Y\!\times\![0,1])}\!=\!T_Y.$ Thus  $\Psi(Y\!\times\![0,1])$ is  a hornic neighborhood of $Y,$ by Example \ref{Ex.Fast.Cycles}.v.
  (Note that $X$ is LNE, being a cone.) Therefore $\Psi$ is a hornic link-homotopy.

      \epr

\subsection{The vanishing rate of a (subanalytic) fast cycle  $Y\sset X$}\label{Sec.Preliminaries.Vanishing.Rate}
\bed
 The (inner) vanishing rate is  $\vanrate_X Y:=\sup_{Y'} \{\dist_X(Y,Y')\}.$
 Here the supremum goes over all (subanalytic) non-thin subgerms $Y'\sseteq Y.$
\eed
Thus $\vanrate_X(Y)$ depends only on the Lipschitz type of the pair $(X,Y),$ and not on the choice of the embedding $X\sset \R^N.$
\bex\label{Ex.Vanishing.Rate}
%\leavevmode\vspace{0.4cm}
\bee[\bf i.]
\item
Let $Y\sset X\sset (\R^N,o)$ be a   fast cycle with connected $\Link[Y],$ and suppose    $T_Y$ is a half-line.  Then $\vanrate_X Y\!=\!\diam_{(\R^N,o)} Y>1.$
\bpr
Take an arc $\ga\sset Y.$ Then $\dist_X(Y,\ga)=\diam_X(Y).$ Thus
\[
\vanrate_X Y\ge \diam_X Y \stackrel{\S\ref{Sec.Notations.Basics.vii}}{=} \diam_{(\R^N,o)} Y>1.
\]
 Vice-versa, take a non-thin subgerm $Y'\sset Y$ satisfying:
 \beq
 \dist_X(Y',Y)\ge \vanrate_X(Y)-\ep>1.
 \eeq
  Then $T_{Y'}\sseteq T_Y,$ hence $T_{Y'}$ is also a half-line. Thus $Y'=\{\ga_i\}$ is a finite set of arcs
 and $\dist_X(Y,Y')=\min_i \tord_X(\ga_i,Y).$
  Finally, for any two arcs $\de_1,\de_2\sset Y$ one has: $\tord_X(\de_1,\de_2)\ge \min\{\tord_X(\de_1,\ga_i),\tord_X(\de_2,\ga_i)\},$ for each $i.$
   Therefore
   \beq
   \diam_{(\R^N,o)} Y\ge \dist_X(Y,Y')\ge \vanrate_X(Y)-\ep.
   \eeq

    Altogether: the vanishing rate is the diameter of the germ.
\epr
\item   Hornic link-homotopies preserve fast cycles, but do not preserve  their vanishing rates.

\item
Take a $\be$-horn $X=\{x^2+y^2=z^{2\be},\ z\ge0\}\sset \R^3 $ and its tilted version $X'=\{(x-2z^\be+2z^\de)^2+y^2=z^{2\be},\ z\ge0\} ,$
 for some $1<\be<\de.$ The two horns intersect along two arcs with parametrizations $(t^\be-t^\de,\pm\sqrt{2t^{\be+\de}-t^{2\de}},t).$
  Then $\Link[X\cup X']$ is two circles intersecting at two points, homotopically it is $\vee^3 S^1.$
   The vanishing rates of $X,X',X\cup X'$ are equal $\be.$ But the vanishing rate of the ``inner circle" is $\frac{\be+\de}{2}.$

\item The higher-dimensional generalization is the hypersurface
\beq
X=\{\sum^j_{i=1}x^2_i=(\sum^{n+1}_{i=j+1}x^2_i)^\be\}\sset \R^{n+1}, \quad \text{for} \quad 2\le j\le n \quad \text{and} \quad   \be>1.
\eeq
 Then $T_X=\{x_1=\cdots=x_j=0\}=\R^{n+1-j}.$
 Thus $\dim_\R T_X<\dim_\R X=n.$ It is helpful to re-embed $X$ as follows:
 \beq
 \tX:=\{(x,t)|\ \sum^j_{i=1}x^2_i=t^{2\be}, \quad \sum^{n+1}_{i=j+1}x^2_i=t^2, \ t\ge0\}\sset \R^{n+1}\times\R^1_t.
 \eeq

\noindent Then $\Link[\tX]$ can be realized as $\tX\cap \{t=t_o\}.$ Thus $\Link[X]\!\cong\! S^{j-1}_\be\!\times\! S^{n-j},$ in particular it does not admit any homotopy to a space of smaller dimension.
  Therefore $X$ is a fast cycle.

 \textbf{Claim.} $\vanrate_X X=\be.$

 \noindent \bpr \underline{The part $\vanrate_X X\ge \be.$}
Take the section $Y=\{x_1=\cdots=x_j\}\cap X.$ Then $\dim_\R Y=n+1-j$ and $T_Y=T_X.$  Thus the subgerm $Y\sset X$ is non-thin.

We claim: $\dist_X(Y,X)\ge\be.$ Indeed, take any arc $\ga\sset X,$ its length parametrization has the form
 $(t^\be c_1(t),\dots,t^\be c_j(t),x_{j+1}(t),\dots,x_{n+1}(t)).$ These functions are subanalytic, of orders$\ge1,$ and at least one of $x_{j+1}(t),\dots,x_{n+1}(t)$ is linear in $t.$
 Accordingly take an arc $\ga_Y(t)=(t^\be z(t),\dots,t^\be z(t),x_{j+1}(t),\dots,x_{n+1}(t)).$
  To ensure that $\ga_Y\sset Y $ we impose the condition $j\cdot z(t)^2=(\sum^{n+1}_{j+1} \frac{x_i(t)}{t})^2.$
   Thus  $z(t)$ is a subanalytic  function of order$\ge1.$ Hence
     $\ga_Y(t)$ is a length parametrization up to $O(t^\be)$. And we have:
   \beq
\tord_X(\ga,\ga_Y)=\ord_t d_X^{sub}(\ga_Y(t),\ga(t))\ge\be.
   \eeq
   Therefore $\vanrate_X X\ge \be.$

\underline{The part $\vanrate_X X\le \be.$} Let $Y\sset X$ be any (subanalytic, closed) non-thin subgerm, thus $\dim_\R Y\le \dim_\R T_X.$
 We claim: $\dist_X(Y,X)\le \be.$ Define the weighted-polar coordinates
 \beq
 \bS\times\R_{\ge0}\stackrel{\si}{\to}\R^{n+1}, \quad \text{ by } \quad
 (\us,t)\to (t^\be s_1,\dots,t^\be s_j,s_{j+1},\dots,s_{n+1}).
  \eeq
  Take the strict transform $\tY:=\overline{\si^{-1}(Y\smin o)}.$
   It is a subanalytic, closed subset  of $\dim_\R=n+1-j.$ Therefore $\dim_\R(\tY\cap (\bS\times \{o\})=n-j<\dim_\R(X\cap (\bS\times\{o\}))=n-1.$
    Thus there exists an arc $\ga\sset X$ for which $\tga\cap (\bS\times\{o\})\not\in\tY.$ But then $\dist_X(\ga,Y)\le \dist_{(\R^{n+1},o)}(\ga,Y)\le \be$. Hence, the Claim is proved.

\eee
\eex

\beR In view of Example \ref{Ex.Vanishing.Rate}.i one could try to define the vanishing rate as $\dist_X(Y,T_Y).$ This is not possible for several reasons.
\bee[\bf i.]
\item
$T_Y\sset (\R^N,o)$ does not lie inside $X.$ Moreover, for $Y=X,$ in many cases $X$ cannot be (subanalytically) re-embedded to ensure $T_Y\sset X.$
\\
But the notion of vanishing rate should depend on the Lipschitz type of the pair $Y\sset X$ only, not on the embedding $X\sset (\R^N,o).$
\item
Let $X=(\R^N,o),$ thus $T_Y\sset X.$ Then
$\dist_{(\R^N,o)}(Y,T_Y)$ is not preserved by analytic coordinate changes. E.g. let $Y=\{x^2+y^2=z^{2\be},\ z\ge0\}\sset \R^3,$ with $\be>2.$ Then $dist(Y,T_Y)=\be.$
 Take the coordinate change $\Phi: (x,y)\to (x+\phi_x(x,y,z),y+\phi_y(x,y,z)),$ where $\phi_x,\phi_y$ are generic polynomials of order $2.$ Then $\Phi(T_Y)=T_Y=\R_{\ge0}\hz.$  But $\dist_{(\R^3,o)}(\Phi(Y),T_Y)=2<\be.$
 Indeed, take an arc $\ga\sset Y$ with parametrization $\ga(t)=(t^\be a_x,t^\be a_y,t),$ for some constants $a^2_x+a^2_y=1.$
  Then $\Phi(\ga)\sset \Phi(Y)$ has the parametrization  $(t^\be a_x+t^2\cdot(\cdots),t^\be a_y+t^2\cdot(\cdots),t).$ Therefore
   $\tord_{(\R^3,o)}(\Phi(\ga),T_Y)=2.$

\item  One could hope to have $\vanrate_X (Y)=\sup\ \dist_X(Y,T_Y),$ the supremum taken over all the coordinate choices for which   $T_Y\sset X.$
 Even this fails.
 E.g. fix $\be>1$ and take the union of two tilted horns in $\R^3:$
\beq
Y:=X:=\{(x_1-t^\be)^2+x^2_2=t^{2\be},\ t\ge0\}\ \cup\ \{(x_1-3t^\be)^2+x^2_2=t^{2\be},\ t\ge0\}.
\eeq
 These horns intersect at $o$ only. Here $T_X=\R_{\ge0}\hat{t}\sset X,$ and $X$ is thin and $\Link[X]=S^1\amalg S^1.$ (Thus $X$ is a fast cycle.) But $\dist_X(T_X,X)=1,$ while $\vanrate_X(X)=\be.$

 In this example $\Link[X]$ is non-connected. A case with connected link is obtained, e.g. by adding the coordinate $x_3$ and a $\dim=2$ plane:
 \[
 \R^4_{x_1x_2x_3t}\supset \cX:=(X\times\R^1_{x_3})\cup \{x_1=0=x_3\}.
 \]
 Here  $\Link[\cX]=\Link[X\times \R^1_{x_3}]\cup S^1_{x_1,t}$ is connected.
  And $\dim_\R \cX=3,$ but $T_\cX=\{x_1=0=x_2\}\cup\{x_2=0=x_3\}.$ Thus $\cX$ is thin. And $\vanrate_\cX \cX=\be.$ But $\dist_\cX(\cX,T_\cX)=1.$
 Moreover, $\dist_\cX(\cX,T_\cX)=1$ for all the coordinate choices.

\eee\eeR

\section{Detecting fast cycles on perturbations of weighted-homogeneous germs}\label{Sec.Fast.Cycles.on.perturbed.weighted.homogen.germs}

\subsection{Germs with (perturbed) weighted-homogeneous foliations}\label{Sec.3.1}   Let $\k=\R,\C.$
 We consider  subanalytic germs $X\sset (\k^N,o)$ foliated by (perturbed) weighted-homogeneous arcs. Namely, we assume a family of arcs $\{\ga_\us\}_{\us\in \Link[X]}$ with (uniform) parametrization
 \beq\label{Eq.foliation.w.hom.arcs}
 \ga_\us (t)= (t^{\om_1}s_1+t^{\om_1}g_1(\us,t),\dots,t^{\om_N}s_N+t^{\om_N}g_N(\us,t)),\hspace{1cm}
 0<\om_1\le \cdots\le \om_N .
 \eeq
 Here:
 \bei
 \item The higher-order-terms  are continuous, $ g_i(\us,t)\in C^0(\Link[X]\times [0,\ep]),$   and vanish at $t=0.$
  (Thus $\liml_{t\to0}\ug(\us,t)\!=\!o.$) As usual, $0<\ep\ll1.$
  \item $X=\cup_{\us\in \Link[X]}\ga_\us$ and $\ga_\us\cap \ga_{\us'}=o$ for $\us\neq \us'.$
 \eei
 Take the $t$-order, $\ord_t g_i(\us,t)\!=\!\sup\{c|\ g_i(\us,t)\!=\!O(t^c)\}.$

\bex\label{Ex.foliations}
%\leavevmode\vspace{0.4cm}
\bee[\bf i.]
\item   Any $\C$-analytic branch $C\sset (\C^N,o)$ admits such a foliation. It is defined e.g. via the normalization map, $(\C^1,o)\to C,$
 where on $(\C^1,o)=(\R^2,o)$ one takes just the radial foliation.

\item
Take a weighted-homogeneous  subvariety  $X_o=V(\uf_\up)\sset \k^N,$ of any dimension.
 The defining polynomials are weighted homogeneous, $f_{p_i}(t^\uom\cdot \ux)=t^{p_i}\cdot f_{p_i}(\ux),$ see \S\ref{Sec.Notations.Basics.iii}.
 Take the foliation of $\k^N$ by weighted-homogeneous arcs, $\ga_\us(t):=t^\uom\cdot\us,$ for $\us\in \bS.$ It is compatible with $X_o,$ i.e. either $\ga_\us\sset X_o$ or $\ga_\us\cap X_o=\{o\}.$
 \eee
 \eex
\bex
Even when $X$ is not weighted-homogeneous, it often has an open (dense) subset  admitting a perturbed weighted-homogeneous foliation.
 E.g. let  $X_o=V(\uf_\up)\sset \k^N$ be a weighted-homogeneous complete intersection, with arbitrary singularities, possibly non-reduced.
 Perturb $X_o$ by (analytic) higher order terms, $X=V(\uf_{\up}+ \uf_{>\up})\sset (\k^N,o).$ One would like to perturb the weighted-homogeneous foliation $\ga_\us$
  into a foliation $\tga_{\us}$ compatibly with $X,$ as above. (Namely,
$X=\cup \widetilde{\gamma}_{\us}
\quad\text{and}\quad
\widetilde{\gamma}_{\us}\cap
\widetilde{\gamma}_{\us'}=\{o\}$
for $\us,\us'\in\Link[X]$ with $\us\neq\us'$.)

 This is not always possible, even if $X_o$ is an ICIS. See Example 2.1 of \cite{Kerner-Mendes.Defs}. However, in that paper we identify the ``obstruction locus" $\Si\sset X_o,$ outside of which such a perturbation  $\ga_\us\rightsquigarrow \tga_{\us}$ does exist. In more detail, refine the weights splitting of \eqref{Eq.weights.split} into
 \beq\label{Eq.weights.split.refined}
 \om_1=\cdots=\om_{r_1 }<\om_{r_1 +1}=\cdots= \om_{r_2}<\om_{r_2+1}=\cdots\le \om_N, \ \text{ for some } \  1\le r_1< r_2<\cdots.
 \eeq
Then the obstruction locus is the union of singular loci of the plane-sections:
\beq
\Si\!: =\!\Sing[X_o]\ \cup\   \Sing[X_o\!\cap\! V(x_1,\dots,x_{r_1})]\  \cup \  \Sing[X_o\!\cap\! V(x_1,\dots,x_{r_2})] \ \cup \dots.
\eeq
\eex
 \noindent
\bel   \cite[Lemma 2.3]{Kerner-Mendes.Defs}  For any neighborhood $\Link[\Si]\sset \cU_\bS\sset \bS$   there exists a hornic neighborhood $\Si\sset \cU(\Si)\sset(\k^N,o)$ with
  $\Link[\cU(\Si)]\sseteq\cU_\bS,$ that satisfies: the germ $(\k^N\smin \cU(\Si),o)$ is foliated by arcs $\tga_{ \us},$ compatibly with $X .$
\eel

Here  $\tga_{ \us}\!=\!t^\uom \cdot \us+  t^\De\cdot t^\uom\cdot \uh(t,\us ),$ with
 $\De=\min_i\{ \ord_t f_{>p_i}(t^\uom \cdot \us)-\ord_t f_{p_i}(t^\uom \cdot \us)\}$ and $\uh\in C^\om([0,t_o)\times (\bS\smin \cU_\bS) ).$

If $\Si=\{o\}$ then $\Link[\Si]=\empty,$ thus the whole $(\k^N,o)$ is foliated by such arcs.
\noindent Some germs admit many different foliations, and the weights are not-uniquely determined,  e.g. one can choose $\om_1\!=\!\om_2$ or $\om_1\!<\!\om_2.$
 Thus we make a choice of the weights.  (See remark \ref{Rem.After.Theorem.4.7}.)

\subsection{}\label{Sec.3.2}
Having such a (perturbed) weighted-homogeneous foliation $\{\ga_\us\}$ on $X,$ we define the semi-group action $\tau\circ \ga_\us(t):=\ga_\us(\tau\cdot t),$
 for $\tau\in (0,1].$
If the foliation is global, e.g. when $X\sset \k^N$ is a weighted-homogeneous subvariety,
 then we get the standard (multiplicative) group-action $\R_{>0}\circlearrowright X.$
 But in general the foliation is only local, thus we get only the action  $(0,1]\circlearrowright X.$

 For a closed subset $Z\sset X,$ or $Z\sset \Link[X],$   we take its orbit under this semi-group action.
  Abusing the notation we denote this orbit by $\R_{>0}Z,$ and take its closure $\R_{\ge0}Z:=\overline{\R_{>0}Z}.$
 This abuse of notation  causes no harm, as we always consider only the germ of $\R_{\ge0}Z$ at $o\in \k^N.$

\subsection{The general criterion for fast cycles}\label{Sec.R-action.IMC.restriction.on.weights.Basic}
 Take a (closed) subanalytic germ $X\sset (\k^N,o).$  Suppose a subanalytic subgerm $Y\sseteq X $ (not necessarily closed or open or dense) admits the foliation of \S\ref{Sec.3.1}.  For any compact subset $Z\sset \Link[Y],$ we get the subgerm $\R_{\ge0}Z \sset Y,$ as in \S\ref{Sec.3.2}.
 Split the weights as in \eqref{Eq.weights.split}, so that
  $\om_1\!=\!\cdots\!=\!\om_r\!<\!\om_{r+1}.$    Take the projection $\pi: \k^N\!\to\! \k^r_{x_1\dots x_r}$ .
\bel\label{Thm.IMC.restriction.weights.general}
    Suppose  for some $0<t_o\ll1$  there exists a cycle $Z\sset Y\cap V(x_1-t_o)$ satisfying:
    \bei
        \item   $\dim_\R Z>\dim_\R \pi Z.$
    \item $Z$ admits no  deformation inside $Y\cap \pi^{-1}\pi(Z)$  to a set of smaller dimension.
    \item
     $Y$ contains the germ of every hornic neighborhood of $\R_{\ge0}Z$ in $X.$
    \eei

   Then the germ $\R_{\ge0}Z\sset X$ is a fast cycle with the following properties.

   \bee%[\!\!\bf 1.\!]
   \item    Its tangent cone lies in the half-plane $\R_{\ge0}\hx_1\times \k^{r-1}_{x_2\dots x_r}. $ \footnote{we denote by $\hx_1$ the unit vector in the first coordinate  direction of $\k^N$}
\\
If $\k=\R$ then $\dim_\R T_{\R_{\ge0}Z}\le r.$

\item The   vanishing rate of this fast cycle satisfies: $\vanrate_X(\R_{\ge0}Z)\ge \frac{\om_{r+1} }{\om_1}>1.$

\item In particular $X$ cannot be IMC.  (By lemma \ref{Thm.IMC.germs.have.no.fast.cycles})
\eee
\eel
\bpr The proof goes in several step.
\bee[Step 1:]
\item  we prove that $\R_{\ge0}Z$ is a thin set with the   tangent cone inside $\R_{\ge0}\hx_1\times \k^{r-1}_{x_2\dots x_r}. $
\item  we prove that  $\R_{\ge0}Z$ does not admit  hornic link-homotopy to a germ of smaller dimension. Hence  $\R_{\ge0}Z$ is a fast cycle.

\noindent It remains to bound the vanishing rate, proving part 2 of the lemma.
\item  we do a warmup case, assuming that the projection of the punctured germ $\R_{>0}Z\to \pi(\R_{>0}Z)\sset \k^r$ is a $C^0$-trivial fibration.
 Then the bound on the vanishing rate is just a direct computation.
\item in reality this projection is never  a $C^0$-trivial fibration. But it is ``stratified-$C^0$-trivial".
 And we verify the bound on the vanishing rate by going over all the strata.
\eee
\medskip

  Rescale the coordinate $x_1 $ and $t$ to  assume:  $Z\sseteq Y\cap V(x_1-1).$   We can assume: $Z$ is compact, connected, subanalytic
   (\S\ref{Sec.Notations.Basics.viii}), of pure dimension.
 Finally, we can assume (by rescaling $x_2\dots x_N$):
 \beq\label{Eq.fast.cycle.location}
 Z\sset X\cap V(x_1-1)\cap \{\|(x_2,\dots, x_N)\|<\ep\ll1\}.
\eeq
\bee[\hspace{-0.6cm}\bf Step 1.]
\item
   The assumed foliation of  $Y$  gives the closed sub-analytic  germ $\R_{\ge0}Z:=\overline{\R_{>0}Z}\sseteq Y\sseteq X.$

Take the set $Box:=\{|x_1|\le1\}\times \{\|(x_2,\dots, x_N)\|\le 1\}.$ For $\k=\R$ this is a cylinder, for $\k=\C$ this is a polydisc.
   Realize $\Link[X]$ inside the boundary of this set,
   \beq
   \Link[X]\cong \di Box\cap X=[X\cap V(|x_1|-1)]\cup [X\cap V(\|(x_2,\dots,x_N)\|-1)].
   \eeq
   Thus $\Link[\R_{\ge0}Z]=\R_{\ge0}Z\cap V(x_1-1)=Z,$ by \eqref{Eq.fast.cycle.location}.    (Both for $\k=\R$ and for $\k=\C$.)

 For the projection $\pi:\k^N\to \k^r_{x_1\dots x_r}$
    we claim: $T_{\R_{\ge0}Z}= T_{\pi(\R_{\ge0} Z)}.$
    \\(In particular
    $T_{\R_{\ge0}Z}\sseteq  \R_{\ge0}\hx_1\times \k^{r-1}_{x_2\dots x_r}. $
       Thus any arc on $\R_{\ge0}Z$ is tangent to this half-plane.)
\\Indeed,  take a (subanalytic) arc $\ga\sset \R_{\ge0}Z$ and its image $\pi\ga\sset \pi\R_{\ge0}Z.$ Each point of $\ga$ belongs to some arc of the foliation. Thus the parametrization of $\ga$ comes from \eqref{Eq.foliation.w.hom.arcs}:
 \beq
 \ga(t)=[t^{\om_1},\ t^{\om_2}(z_2(t)+g_2(\uz(t),t)),\ \dots\ ,t^{\om_N}(z_N(t)+g_N(\uz(t),t))].
 \eeq
Here the function $\uz(t)$ can be discontinuous, but it is bounded (because $Z$ is compact). In particular, $\lim_{t\to0}g_i(\uz(t),t)=0$ for all $i.$

 Moreover there exists $\uz(o):=\lim_{t\to o}\uz(t).$
  Indeed, pass to the weighted-homogeneous polar coordinates,
    $\bS\times\R_{\ge0}\stackrel{\phi}{\to}(\k^N,o),$ by $(\us,t)\to(s_1 t^{\om_1},\dots,s_N t^{\om_N}).$
     The strict transform of $\ga$ is a subanalytic arc $\phi^*\ga.$   And  $\lim_{t\to o}\uz(t)=\phi(\lim_{t\to o}\phi^*\ga(t))$ exists.

 Therefore the tangent line to $\ga$ is spanned by the vector $(1,z_2(o),\dots,z_{r }(o),0,\dots,0)\in \R_{\ge0}\hx_1\times \k^{r-1}_{x_2\dots x_r}.$ And the same vector spans the tangent line to $\pi(\ga).$
\\We have proved: $T_{\R_{\ge0}Z}\sseteq T_{\R_{\ge0}\pi(Z)}.$ The inclusion $T_{\R_{\ge0}Z}\supseteq T_{\R_{\ge0}\pi(Z)}$ is proved similarly.

\medskip

  The subset $\pi \R_{\ge0}Z\sset \R_{\ge0}\hx_1\times \k^{r-1}_{x_2\dots x_r}$ is closed, connected, subanalytic.
  It is covered by the ``almost radial" foliation, $\ga_\uz(t)=t\cdot (1,z_2+g_2(\uz,t),\dots,z_r+g_r(\uz,t)),$ with $\liml_{t\to o}g_i(\uz,t)=0.$
  Therefore:
  \beq
  \dim_\R T_{\R_{\ge0}Z}=\dim_\R T_{\pi(\R_{\ge0}Z)}=  \dim_\R \pi (\R_{\ge0} Z)=\dim_\R \pi (Z)+1.
  \eeq
By the assumption we have  $\dim_\R \pi (Z)\< \dim_\R Z.$
   Thus $\R_{\ge0}Z$ is a thin subset.

   \medskip

   \item
  We claim: $\R_{\ge0}Z$ does not admit any  hornic link-homotopy (inside $X$) to a germ of a smaller dimension.

It is enough to prove: any hornic neighborhood $\cU(\R_{\ge 0}Z)\sset X$ admits a hornic deformation-retraction to a foliated subgerm
  $\tcU\sset Y$  that satisfies: $\pi(\tcU\cap V(x_1-t_o))\sseteq \pi (Z).$ Indeed, assuming this, any hornic link-homotopy of
   $\R_{\ge0}Z$ can be done inside such $\tcU.$ But this restricts to a homotopy of $Z$ inside the set $\pi^{-1}\pi(Z)\cap Y.$
    And such a homotopy cannot decrease $\dim_\R Z,$ by our assumption.

  Thus we fix a hornic neighborhood $\cU(\R_{\ge 0}Z)\sset X.$
 By our assumption: $\cU(\R_{\ge 0}Z)\sset Y.$
 Now we repeat the proof of Lemma \ref{Thm.IMC.germs.have.no.fast.cycles} for our weighted homogeneous case. Take a contracting flow  $\Phi^\cU_\tau\circlearrowright \cU(\R_{\ge 0}Z) ,$ and the ``perturbed weighted-homogeneous flow"  of \S\ref{Sec.3.2},
  $\Phi^{Y}_\tau\circlearrowright Y .$ Their composition, $\Psi_\tau:=\Phi^{Y}_\frac{1}{\tau}\circ \Phi^{\cU}_\tau,$ is a link-homotopy of
   $\cU(\R_{\ge 0}Z)$ inside $Y,$ defined for $\tau>0.$ It extends to $\tau=0$ continuously. (As in the proof of Lemma \ref{Thm.IMC.germs.have.no.fast.cycles}, we pass to the weighted-homogeneous coordinates.)

 Finally, the link-homotopy $\Psi_\tau$ is hornic, i.e. it preserves the tangent cone.
   Therefore for any point $p\in \cU(\R_{\ge 0}Z)\cap V(x_1-t_o)$ we get: $lim_{\tau\to 0}\pi\Psi_\tau(p)\in T_{\pi(\R_{\ge0}Z)}.$

 \medskip

\hspace{-1cm}Together with Step 1 we get:
  $ \R_{\ge0}Z\sset X$ is a fast cycle of $\dim_\R \R_{\ge 0}Z=\dim_\R Z+1.$

\hspace{-1cm}It remains to prove: $\vanrate_X \R_{\ge0}Z\ge \frac{\om_{r+1}}{\om_1}.$ It is enough to construct a (subanalytic,

\hspace{-1cm}compact)  subset $\tZ\sset Z$ such that $\R_{\ge0}\tZ$ is non-thin and satisfies:
\beq\label{Eq.inside.proof.conditions.to.ensure}
T_{\R_{\ge0}\tZ}=T_{\R_{\ge0}Z},\hspace{1.2cm}
\dim_\R T_{\R_{\ge0}\tZ}=\dim_\R \R_{\ge0}\tZ, \hspace{1.2cm}
\dist_X(\R_{\ge0}\tZ,\R_{\ge0}Z)\ge \frac{\om_{r+1}}{\om_1}.
\eeq

\item
  As a warmup we assume that the projection of punctured germ $\R_{>0}Z\!\to\! \pi(\R_{>0}Z)\sset \k^r$ is a $C^0$-trivial fibration.
 That is, we assume that there exists a (subanalytic) homeomorphism  $\Psi\!: \R_{>0}Z\!\isom{}\!\pi(\R_{>0}Z)\!\times \! F,$
 with $\pi$ acting as a projection onto the first factor. (Here $\Psi$ is not necessarily bi-Lipschitz.) The fibre $F$ has a finite number of connected components (by subanalyticity).
  Fix a point in each connected component of the fibre, $\{y_\bullet\}\sset F.$
  This defines the ``multi-section of $\pi$",
  \beq
  \R_{\ge0}\tZ:=\overline{\Psi^{-1}\big(\pi(\R_{>0}Z)\times \{y_\bullet\}\big)}\sset \R_{\ge0}Z.
  \eeq
\noindent\parbox{11.5cm}
  {We claim: the set $\R_{\ge0}\tZ$ satisfies
   the conditions of \eqref{Eq.inside.proof.conditions.to.ensure}. Indeed, the closedness, subanalyticity and the first two equalities of \eqref{Eq.inside.proof.conditions.to.ensure} are obvious.
 To bound $\dist_X(\R_{\ge0}\tZ,\R_{\ge0}Z)$ we take any arc $\ga\sset \R_{\ge0}Z.$ Then $\Psi^{-1}(\pi\ga\times\{y_\bullet\})\sset  \R_{\ge0}\tZ$ is  a finite set of arcs.  Exactly one of them (denote it $\tga$) lies in the same connected component of $\pi^{-1}\pi(\ga)$ as $\ga.$
 }
  \begin{picture}(0,0)(-15,-33)
  \color{gray}\thicklines
 \curve(1,0.2,10,2,20,4,30,7,40,11) \curve(1,-0.2,10,-2,20,-4,30,-7,40,-11)

  \curve(1,0.5,10,3,20,9,30,19,40,31)  \curve(1,0.5,10,5,20,14,30,30,40,50)

   \curve(1,-0.5,10,-3,20,-9,30,-19,40,-31)  \curve(1,-0.5,10,-5,20,-14,30,-30,40,-50)
\thinlines
\multiput(20,4)(5,1){4}{\line(1,-1){10}}
\curve(35,26,40,46)\curve(31,22,36,39)\curve(27,17,31,30)
\curve(35,-26,40,-46)\curve(31,-22,36,-39)\curve(27,-17,31,-30)

\color{black}
\put(-3,13){\vector(0,-1){30}}  \put(-12,-2){$\pi$}  \put(0,20){$\scriptstyle{\R_{>0}Z}$}

\curve(0,-70,40,-70)
\put(50,-70){$\scriptstyle{\pi(\R_{>0}Z)}$}

\color{red}\thicklines \curve(1,0,10,0,45,0)  \put(47,-4){$\scriptstyle{\pmb\ga}$}

\thinlines\color{blue}
\curve(1,0,10,1,20,2,45,6)
\curve(1,0,10,3,20,11,45,44)
\curve(1,0,10,-3,20,-11,45,-44)     \put(47,30){$\scriptstyle{\pmb{\psi^{-1}(\pi\ga\times\{y_\bullet\})}}$}
\put(47,4){$\scriptstyle{\pmb\tga}$}
\end{picture}

\medskip

\noindent\mbox{By the construction $\tord_{(\k^N,o)}(\ga,\tga)\!\ge\! \frac{\om_{r+1}}{\om_1}. $
 Moreover, we claim: $\tord_{\R_{\ge0}Z}(\ga,\tga)\!\ge \!\frac{\om_{r+1}}{\om_1}. $}
   \\Indeed, fix parametrizations $\ga(t),\tga(t)$  such that $\pi (\ga(t))=\pi (\tga(t)).$ Then (for each $0<t_o\ll1$) the points $\ga(t_o),\tga(t_o)$ are connected by a path
   \beq
   [\ga(t_o)\stackrel{\eta_{t_o}}{\rightsquigarrow}\tga(t_o)]\sset \pi^{-1}(\pi \ga(t_o))\cap \R_{\ge0}Z.
   \eeq
  The parametrization of this path is
  \beq
  \eta_{t_o}(s)=(t_o\cdot z_1(t_o),\dots,t_o\cdot z_r(t_o),t^{\frac{\om_{r+1}}{\om_1}}_o\cdot z_{r+1}(t_o,s),\dots, t^{\frac{\om_N}{\om_1}}_o\cdot z_N(t_o,s)).
   \eeq
   Here  $z_1(t_o) \dots z_r(t_o)$ are constants  (for a fixed $t_o$), while $z_{r+1}(t_o,s),\dots,z_N(t_o,s)\!\in\! C^0([0,1]).$
   Moreover, by subanalyticity of $\uz,$ we can assume: the function $\eta_{t_o}$ (and   $z_{r+1}(t_o,s),\dots,z_N(t_o,s)$) is piecewise $C^1$ and
 with bounded derivatives.
   Invoke the action of \S\ref{Sec.3.2} to get the functions
\beq
z_1(t ),\dots,z_r(t ),z_{r+1}(t ,s),\dots,z_N(t ,s)\in C^0([0,t_o]\times[0,1]).
\eeq
These functions are uniformly continuous.
    Therefore
    \beq
  \quad   length[\eta_t]\=\int^1_0 \sqrt{\sum^N_{i=r+1}t^{\frac{2\om_j}{\om_1}}\cdot
   \Big(\frac{\di z_j(t,s)}{\di s}\Big)^2}ds \le  C\cdot t^{\frac{\om_{r+1}}{\om_1}}.
    \eeq

  Altogether $\dist_X(\R_{\ge0}\tZ,\R_{\ge0}Z)\ge \frac{\om_{r+1}}{\om_1}.$

\item
 In reality the projection of the punctured germ,  $\R_{>0}Z\to \pi(\R_{>0}Z),$ is never a $C^0$-trivial fibration.
 But it is subanalytic and proper. Hence the triviality holds outside of the branching locus. This branching locus is subanalytic, of $\codim_\R\ge1.$
  Therefore we take a (subanalytic) subset $W\sseteq \pi(\R_{>0}Z)$ satisfying:
\bei
\item $W$ is open, dense inside $\pi(\R_{>0}Z),$ and has a finite number of connected components, $W=\amalg_j W_j$
\item $\pi:\R_{>0}Z\to \pi(\R_{>0}Z)$ is a trivial proper fibration over each $W_j$, i.e. ${\psi_j}: \pi^{-1}W_j\cap \R_{>0}Z\isom W_j\times F_j,$
 with $\pi$ acting as the projection onto $W_j.$
\eei
Treat each $W_j$ as in the warmup, fixing a point on each connected component of the fibre, $\{y_\bullet\}\sset F_j.$ Take the corresponding ``multi-section" $\psi^{-1}(W_j\times\{y_\bullet\})\sset \R_{>0}Z.$ One gets the subset $\R_{\ge0}Z_0:=\amalg_j \overline{\psi^{-1}(W_j\times\{y_\bullet\})}\sset \R_{>0}Z$ satisfying:
\bei
\item
$T_{\R_{\ge0}Z_0}=T_{\R_{\ge0}W};$ \quad  $\dim_\R T_{\R_{\ge0}Z_0}=\dim_\R  \R_{\ge0}Z_0 ;$ %(thus $\R_{\ge0}Z_0$ is not thin);

\item
$\dist_X(\R_{\ge0}Z_0,\R_{\ge0} Z\cap \pi^{-1}W) \ge \frac{\om_{r+1}}{\om_1}.$
\eei

Restrict $\pi$ to the complement $(\R_{\ge0} Z)\smin \pi^{-1}W.$ This set is subanalytic, closed, with $\dim_\R (\R_{\ge0} Z\smin \pi^{-1}W)< \dim_\R \R_{\ge0}  Z.$
 Repeat this construction for $(\R_{\ge0}  Z\smin \pi^{-1}W)$ and pass to the further complements.
  The process stops after at most $\dim_\R   Z$-steps. We get the subset $\tZ:=\cup Z_i\sset Z $  that satisfies:
 \bei
 \item     $\pi(\tZ)=\pi(Z)$  and $\pi(\R_{\ge0}\tZ)=\pi(\R_{\ge0}Z).$
\item
$T_{\R_{\ge0}\tZ }=T_{\R_{\ge0}Z};$ \quad  $\dim_\R T_{\R_{\ge0}\tZ}=\dim_\R  \R_{\ge0}\tZ ,$ thus   $\R_{\ge0}\tZ$ is not thin.
\item
$\dist_X(\R_{\ge0}\tZ,\R_{\ge0}Z)\ge \frac{\om_{r+1}}{\om_1}.$
\eei
This proves the statement. \epr
\eee

\subsection{Fast cycles on perturbations of complex weighted-homogeneous complete intersections}\label{Sec.R-action.IMC.restriction.on.weights.Full} In this subsection $\k=\C.$
Let $X_o :=V(\uf_\up)\sset (\C^N,o)$ be   a   (complex-analytic) weighted-homogeneous complete intersection  of dimension $n\ge2$, with
 arbitrary singularity, possibly non-reduced.
 We assume   $(\uf_\up)\sseteq (\ux)^2.$
   Take a perturbation by terms of (strictly) higher orders, $X\!:=\!V(\uf_\up+\uf_{>\up})\sset (\C^N,o).$ We assume: $X$ is reduced, i.e.
   the ideal $(\uf_\up+\uf_{>\up})\sset \C\{x\}$ is radical.
    Take small representatives of these germs.

\subsubsection{The surface case, $n=2$}
\bprop\label{Thm.Fast.Cycles.via.weights.Surfaces}
Suppose for each $0<t_o\ll1$ the complex affine curve  $X\cap V(x_1-t_o)$ contains a connected non-contractible component.
 If $X$ is not weighted-homogeneous, then suppose (in addition) that this component is locally irreducible at each point. If  $\om_1<\om_2$ then $X$ has a fast loop of vanishing rate$\ge\frac{\om_2}{\om_1}.$
\eprop
\bpr
     The assumption implies a loop $Z,$ non-contractible  inside  $X\cap V(x_1-t_o).$
    \bei
    \item   If $X$ is weighted-homogeneous, then we get the germ $\R_{\ge0}Z\sset X.$

\item In the non-weighted-homogeneous case we have the  deformed foliation on $X,$ off the obstruction locus $\Si,$ by example \ref{Sec.3.1}.iii.
  Note that $\Si\sset \Sing(X_o)\cup [X_o\cap V(x_1)].$
             By a small deformation of $Z$ (inside the  $\C$-analytic curve $X\cap V(x_1-t_o),$ which is locally irreducible at each point)  we can assume: $Z$ does not intersect
              the subset $X\cap \Sing(X_o)\cap V(x_1-t_o).$  Indeed, as  $X$ is reduced of $\dim_\C=2$, this intersection is a finite set of points.

    Note that all the other components of $\Si$  do not intersect $V(x_1-t_o).$
Therefore, pushing  $Z$ faraway from  $X\cap \Sing(X_o)\cap V(x_1-t_o),$  we can assume: $Z$ lies inside the part of $X$ that is foliated as in Example \ref{Sec.3.1}.iii.
\eei

In both cases we get the subgerm $\R_{\ge0}Z\sset X.$

Finally, to invoke lemma \ref{Thm.IMC.restriction.weights.general} we observe:   $r=1$ and  $\dim_\R Z> \dim_\R \pi Z=0.$
 Thus we get the fast loop $\R_{\ge0}Z,$ of the claimed vanishing rate.
\epr

\subsubsection{Arbitrary dimension, $n\ge2$}  Suppose  $X\!\cap\! V(x_1)$ is a reduced complete intersection of dimension $(n-1).$
  In the general case    various pathologies can occur, see Remark \ref{Rem.After.Theorem.4.7}. We make two assumptions:

    \bee[\quad \bf i.]%\hspace{-0.3cm}
    \item \
  [For the case when $X$ or $X\cap V(x_1)$ have non-isolated singularities]
\\
Suppose  $\Sing(X)\sset V(x_1), $ and for each $0\!<\!t_o\!\ll\!1$ the section  $X\!\cap\! V(x_1-t_o)$ is the (smooth) Milnor fibre of the germ $X\cap V(x_1).$

     \item \ [For the case $X\neq X_o,$ i.e. the perturbation $X_o\rightsquigarrow X$ is non-trivial]  \
     \\Suppose   $\dim_\C[X\cap \Sing(X_o)\smin V(x_1) ]<\frac{n+1}{2}.$
     \eee
 Observe that the assumption i. is irrelevant when both $X$ and $X\cap V(x_1)$ are ICIS.
 Similarly, the assumption ii. is irrelevant when $X$ is weighted-homogeneous. (In this case we   take  $X_o=X,$ i.e. the perturbation becomes trivial.)

\medskip

Denote $l\!:=\!n-\dim_\C \Sing[X\!\cap\! V(x_1)].$ Thus $2\le l\le n.$
 \bthe\label{Thm.IMC.criterion.via.weights}     Suppose $\om_1\!<\!\om_l, $ thus $1\!\le\! r \!<\!l\!\le\! n.$

Then $X$ has a fast cycle $\R_{\ge0}Z,$  of   dimension $l \! \le\! \dim_\R \R_{\ge0}Z \!\le\! n ,$ and with the properties:
\bee
\item Its tangent cone   is of dimension $\dim_\R T_{\R_{\ge0}Z}\!\le\! r $ and  lies in the half-plane $\R_{\ge0}\hx_1\!\times\! \C^{r-1}_{x_2\dots x_r}.$

  \item     The (hornic) homotopy type of $Z$ is that of the Milnor fibre of the germ $X\cap V(x_1).$
\item  The vanishing rate is bounded, $\vanrate_X(\R_{\ge0}Z)\ge\frac{\om_{r+1}}{\om_1}.$

\item Suppose $X\!\cap\! V(x_1)$ has an isolated singularity, thus $l\=n\!>\!r.$
 If all the fibres of the projection $X\!\to \!(\C^{r+1}_{x_1\dots x_{r+1}},o)$ are of $\dim_\C\!\le\! n-r-1,$ then the vanishing rate is exactly $\frac{\om_{r+1}}{\om_1}.$
 \eee
  \ethe

\bcor
Suppose $X$ satisfies the assumptions i., ii. above. If $X$ is an inner metric cone, then $\om_1=\cdots=\om_l$.
\ecor
\bpr The proof of the Theorem goes in several steps.
 \bee[Step 1:]
\item the space $X\cap V(x_1-t_o)$ is $(l-2)$-connected, and its only non-vanishing  homologies can occur in the (real) dimension-range $[l-1,n-1].$

\item
    we get a non-contractible cycle $Z\sset X\cap V(x_1-t_o),$ lying in the foliated part, and of dimension $l-1\le \dim_\R Z\le n-1.$
    This gives the needed germ $\R_{\ge0}Z.$
\item we invoke Lemma \ref{Thm.IMC.restriction.weights.general} to conclude that  $\R_{\ge0}Z$  is a fast cycle with properties 1-3.
\item we compute the vanishing rate, property 4.
\eee

 \bee[\hspace{-0.5cm}\bf Step 1.]
 \item
 Take the reduced complete intersection $ X\cap V(x_1)\sset (\C^N,o).$
 We claim: $X\cap V(x_1-t_o)$ is  its  smooth  Milnor fibre.
  \bei
  \item If  $X\cap V(x_1)$ has a non-isolated singularity, then this is the initial assumption.
  \item    Otherwise
  $X\cap V(x_1)$ is an ICIS,  and $X\cap V(x_1-t_o)$ is smooth, by Sard's theorem. And thus it is the Milnor fibre of $X\cap V(x_1)$, because the complement of the discriminant in the miniversal deformation is connected.

   The connectedness holds because the miniversal deformation is smooth, and the discriminant is of $\codim_\C=1,$ see e.g. \cite[pg. 24]{AGLV.II}.
\eei

\medskip

\noindent  By \S\ref{Sec.Notations.Milnor.Fibre.Homology} the Milnor fibre  $X\cap V(x_1-t_o)$   has the following properties.
  \bei
  \item It is $(l-2)$-connected, in particular connected.
  \item It admits deformation-retraction to a subset of $\dim_\R\le (n-1).$
 Therefore   the only non-vanishing (reduced) homologies of $X\cap V(x_1-t_o)$ can occur in the (real) dimension-range $[l-1,n-1].$
  \item  It  is non-contractible, and either $\pi_1( X\cap V(x_1-t_o))\neq 0$ or at least one of the homologies
   $H_{l-1}(X\cap V(x_1-t_o),\Z),$ \dots,   $H_{n-1}(X\cap V(x_1-t_o),\Z)$  does not vanish.
  \eei

\item After triangulating $X \cap V(x_1 - t_o)$, choose a compact connected cycle
$Z \subset X \cap V(x_1 - t_o)$ of the minimal possible dimension that represents
a non-zero class either in $\pi_1\big(X \cap V(x_1 - t_o)\big)$ or in
$H_{\dim_{\mathbb R} Z}\big(X \cap V(x_1 - t_o), \mathbb Z\big)$.
We may assume that $Z$ is subanalytic, see \S2.1.viii.
Thus $l - 1 \le \dim_{\mathbb R} Z \le n - 1$.

\medskip

\noindent  We claim:  $Z$ can be chosen to lie inside the foliated part of $X.$
\bei
\item If $X$  is weighted-homogeneous, i.e.  $X=X_o,$ then the whole $X$ is foliated.
\item In the non weighted-homogeneous case  we have the  deformed foliation on  $X,$ off the obstruction locus  $\Si,$ by example \ref{Ex.foliations}.iii.
  Note that $\Si\sset \Sing(X_o)\cup [X_o\cap V(x_1)],$ thus   the part $X\smin \cU[\Sing(X_o)\cup V(x_1)]$
 is foliated.    Here $\cU[\Sing(X_o)\cup V(x_1)]$ is a small hornic neighborhood.

We claim: $Z$ can be deformed outside of $\cU[\Sing(X_o)\cup V(x_1)].$ It is enough to prove: $Z$ can be deformed off the locus
   $X\cap \Sing(X_o)\cap V(x_1-t_o).$ (Then we can push $Z$ away from this locus, inside the manifold $X\cap V(x_1-t_o)$.)
   And this latter claim follows by the general position theorem for definable spaces, e.g. \cite{Loi}, Theorem 2.
    Note the bound on dimensions
    \beq
    \dim_\R Z+\dim_\R[X\cap \Sing(X_o)\cap V(x_1-t_o)]<2(n-1)= \dim_\R(X\cap V(x_1-t_o)).
    \eeq

\eei

\medskip

\item For the projection $X\stackrel{\pi}{\to}\C^r_{x_1\dots x_r}$ we can assume  $\dim_\R\pi Z\le r-1.$ Indeed, we have the morphism of   complex manifolds, $X\cap V(x_1-t_o)\stackrel{\pi}{\to}\pi V(x_1-t_o)\approx\C^{r-1}_{x_2..x_r}.$ The manifold $X\cap V(x_1-t_o)$ is analytic, not necessarily algebraic. But $X\cap V(x_1-t_o)$  can be approximated (near the compact subset $Z$) by a complex-algebraic manifold.
      And then we invoke  the Andreotti-Frankel-Hamm theorem
  for {\em algebraic} morphisms,   \cite{Kerner}. It reads: there exists a deformation-retraction $\{\Phi_t\}_t$ of $\pi$ that acts fibrewise,   making the diagram to commute
 \[\bM
(X\cap V(x_1-t_o))\times[0,1]&\stackrel{\{\Phi_t\}}{\to}&X\cap V(x_1-t_o)
\\(\pi,Id)\downarrow&& \downarrow \pi\\
\C^{r-1}\times[0,1]&\stackrel{\{\pi\circ\Phi_t\}}{\to}&\C^{r-1}
\eM\ ,
\]
 and satisfying: $\dim_\R \pi\circ \Phi_1(\C^{r-1})\le r-1.$    Thus we   replace the original cycle $Z$ by $\Phi_1(Z).$

\medskip

Altogether, $Z\sset X$ lies in the foliated part of $X,$ admits no deformation-retraction to a subset of smaller dimension, and $\dim_\R Z\ge l-1> \dim_\R \pi(Z)=r-1.$ Therefore we invoke Lemma \ref{Thm.IMC.restriction.weights.general} for the cycle $Z$. This proves the statements 1,2,3.

  \item (The exact value of the vanishing rate, statement 4.) In this case $l=n,$ thus $\dim_\R \R_{\ge0}Z=n$ and $0\neq [Z]\in H_{n-1}(X\cap V(x_1-t_o),\Z).$

       We prove: if     $\vanrate_X(\R_{\ge0}Z)>\frac{\om_{r+1}}{\om_1},$ then
      $0= [Z]\in H_{n-1}(X\cap V(x_1-t_o),\Z),$ which is a contradiction. This goes in several steps.

\bee[\bf i.]
\item Suppose        $\vanrate_X(\R_{\ge0}Z)>\frac{\om_{r+1}}{\om_1}.$
   Then there exists a (subanalytic, closed) non-thin subgerm
   $\R_{\ge0}\tZ\sset \R_{\ge0}Z,$ of $\dim_\R (\R_{\ge0}\tZ)\le r,$ that  satisfies: $\dist_X(\R_{\ge0}Z ,\R_{\ge0}\tZ )>\frac{\om_{r+1}}{\om_1}.$
   For the projection $(\C^N,o)\stackrel{\pi_r}{\to} (\C^r_{x_1\dots x_r},o),$  we have:
    \bei
    \item Each fibre of $\R_{\ge0}Z\stackrel{\pi_r}{\!\to\!} (\C^r,o)$ has a finite number of connected components.

    \item   $\dist_{(\C^r,o)}(\pi_r \R_{\ge0}Z,\pi_r \R_{\ge0}\tZ)>\frac{\om_{r+1}}{\om_1}.$
    \eei
Moreover, we can assume that  the fibres of the projection  $ \R_{\ge0}\tZ\to (\C^r,o)$ are    finite.
 (By the same argument as in Step.3 of the proof of Lemma \ref{Thm.IMC.restriction.weights.general}.)

\item   Factorize $\pi_r$ into $X\stackrel{\pi_{r+1}}{\to}(\C^{r+1}_{x_1\dots x_{r+1}},o)\to (\C^r,o). $ Thus the projection
 $\R_{\ge0}Z\to \pi_r\R_{\ge0}Z  $ factors into
   $\R_{\ge0} Z\to\pi_{r+1}\R_{\ge0}Z\to \pi_r\R_{\ge0}Z.$
    We claim: all the fibres of the projection $\pi_{r+1}\R_{\ge0}Z\to \pi_r\R_{\ge0}Z$ are finite.
     Namely, for each point $(s_1,\dots,s_r)\in \pi_r(Z)$ there exists a  finite number of values of $s_{r+1}$ for
      which $(s_1,\dots,s_{r+1})\in \pi_{r+1}(Z).$
   Indeed,   each arc
   \beq
   \ga_s(t)=(t^{\om_1}(s_1 +\cdots),\dots, t^{\om_{r+1}}(s_{r+1} +\cdots), \dots)\sset \R_{\ge0}Z
   \eeq
   is approximated by an arc
 \beq
 \tga_s(t)=(t^{\om_1}(s_1 +\cdots),\dots, t^{\om_{r+1}}(s_{r+1} +\cdots),\dots)\sset \R_{\ge0}\tZ
 \eeq
  with the same constants $s_1,\dots,s_{r+1}.$
 And the set of values of $s_{r+1}$ for $\tga_s$ is finite, as the fibres of  $ \R_{\ge0}\tZ\to (\C^r,o)$ are    finite.

 Therefore $\dim_\R \pi_{r+1}\R_{\ge0}Z= \dim_\R \pi_r\R_{\ge0}Z=\dim_\R T_{\pi_r\R_{\ge0}Z}\le r.$

\item
  Finally, we apply the Andreotti-Frankel-Hamm theorem  for families:

{\em Take a morphism of complex spaces, $X\to B.$ Suppose all its fibres are Stein spaces of $\dim_\C\le n.$ Then this family admits
 a fibrewise deformation-retraction $X\times[0,1]\to Y$ to a subfamily $X\supset Y\to B,$ with all fibres of $\dim_\R\le n.$}

 Indeed, take some stratification making $X_o\to B$ a $C^0$-stratified-trivial morphism, and invoke the classical version,  \cite{Hamm.83},  successively over the strata.

  \medskip

 In our case all the fibres of $X\cap V(x_1-t_o)\to \pi_{r+1}V(x_1-t_o)$ are Stein spaces of $\dim_\C \le n-r-1.$
  Therefore there exists a fibrewise deformation-retraction $X\cap V(x_1-t_o)\stackrel{\Psi}{\rightsquigarrow} Y\cap V(x_1-t_o),$ such that all the fibres of $Y\cap V(x_1-t_o)\to \pi_{r+1}(Y\cap V(x_1-t_o))$ are of $\dim_\R\le n-r-1.$ Therefore
   \beq
\\dim_\R \Psi(Z) \le \underbrace{\dim_\R \pi_{r+1}(Z)}_{\le r-1} + \underbrace{(n-r-1)}_{\text{fibre dimension}} = n-2< \dim(Z)=n-1.
\eeq
    Contradicting the assumption  $0\neq[Z]\in H_{dim(Z)}(X\cap V(x_1-t_o),\Z).$
    \epr
\eee
 \eee

\beR\label{Rem.After.Theorem.4.7}
   In Theorem \ref{Thm.IMC.criterion.via.weights}
 we assume  ``$X\cap V(x_1-t_o)$ is the Milnor fibre",  for the case ``$X\cap V(x_1)$  has a non-isolated singularity". This condition is non-empty.

 E.g. take the 3-fold germ $X=V(x_1x_4+x^d_2+x^d_3)\sset (\C^4,o)$ for $d\ge3.$ Here
 the weights $\om_1,\om_4$ are not uniquely defined.
 \bei
 \item We can put $\om_1<\om_2=\om_3=\frac{1}{d}<\om_4=1-\om_1.$

 Then  $X\cap V(x_1)=V(x_1,x^d_2+x^d_3)\sset (\C^4,o)$ is reduced and $X\cap V(x_1-t_o)$ is smooth (and contractible).
  But $X\cap V(x_1-t_o)$  is not the Milnor fibre of $X\cap V(x_1).$  Thus Theorem \ref{Thm.IMC.criterion.via.weights} does not give any fast cycle.
\item  We can put $\om_2=\om_3=\frac{1}{d}<\om_1<\om_4=1-\om_1.$ Then e.g. $X\cap V(x_2)$ is an isolated hypersurface singularity.
 Then $X\cap V(x_2-t_o)$ is its Milnor fibre, with $H_2(X\cap V(x_2-t_o),\Z)\neq0.$ Thus Theorem \ref{Thm.IMC.criterion.via.weights} ensures a fast cycle.

 \eei
\eeR

\subsection{Examples and applications}
\subsubsection{Semi-weighted-homogeneous ICIS}
 Let $X_o=V(\uf_\up)\sset (\C^N,o)$ be a weighted-homogeneous ICIS  of $\dim_\C=n\ge2.$
  Split the weights as in \eqref{Eq.weights.split}.
 Take a (possibly trivial) perturbation by higher order terms,
  $X:=V(\uf_\up+\uf_{>\up})\sset (\C^N,o).$   Suppose $ X\cap V(x_1)$ is an ICIS of dimension $(n-1).$
   Take the Milnor number $\mu:=\mu(X\cap V(x_1)).$
 \bcor\label{Thm.IMC.criterion.semi.weighted.homogen.germs}
\bee
\item If $\om_1<\om_n$ (i.e.   $r < n$), then $X$ has a fast cycle of (hornic) homotopy type
  $\vee^{\mu}S^{n-1},$  whose tangent cone lies in the half-plane $\R_{\ge0}\hx_1\times \Span_\C(\hx_2,\dots,\hx_{r }),$
  is of $\dim_\R\le r,$
  while the vanishing rate is   $\ge \frac{\om_{r+1}}{\om_1} .$  %%?? Should be equality here.

If all the fibres of the projection $X\to (\C^{r+1}_{x_1\dots x_{r+1}},o)$ are of $\dim_\C\le n-r-1,$ then the vanishing rate is exactly $\frac{\om_{r+1}}{\om_1}.$
\item
  If $X$ is an IMC  then $\om_1=\cdots=\om_n.$
  \eee
  \ecor
  \bpr
   Part 1 follows straight from   theorem \ref{Thm.IMC.criterion.via.weights}.
Part 2 follows by lemma \ref{Thm.IMC.germs.have.no.fast.cycles}.
\epr
 \bex\label{Ex.Examples.of.Fast.Cycles}
 \bee[\bf i.]  %\hspace{-0.5cm}
 \item Part two of this corollary extends  the criterion of \cite{Birbrair-Fernandes.Neumann.08}: ``If a cyclic quotient germ  $\quots{(\C^2,o)}{\mu} $ is IMC,  then its two lowest weights are equal".
 \item
 Consider the Brian\c{c}on-Speder family, $f_\ep(x,y,z)=z^5+x^{15}+xy^7+\ep\cdot zy^6.$
 It is weighted homogeneous with weights $(\frac{1}{15},\frac{2}{15},\frac{3}{15}).$  The Milnor number is $\mu=364.$
  The projection $V(f_\ep)\to\C^2_{xy}$ is finite.

  In this case $V(f_\ep)\cap V(x)=V(z^5+\ep\cdot zy^6,x)$ is reduced for $\ep\neq0.$ Therefore (for $\ep\neq0$) we get: the germ $V(f_\ep)\sset (\C^3,o)$ has a fast cycle of homotopy type $\vee^{364}S^1,$
   whose tangent cone is the half-line $\R_{\ge0}\hx,$ and the vanishing rate  equals $\frac{\om_2}{\om_1}=2.$

 \item (Semi-Brieskorn-Pham singularities, the hypersurface case, $c=1$, $n\ge2.$) Take $f_p(\ux)=\sum   x^{p_i}_i,$ where $p_1=\cdots=p_r> p_{r+1}\ge\cdots  \ge p_{n+1},$  for some $r<n.$
  The hypersurface $V(f_p)\sset (\C^{n+1},o)$ has an isolated singularity.
  Let $f_{>p}(x)$ be of  order$>1$ \wrt the weights $\frac{1}{p_1},\dots, \frac{1}{p_{n+1}}.$
     Then $V(f_p+f_{>p})$ has a fast cycle of homotopy type $\vee^\mu S^{n-1},$ where $\mu=\prod^{n+1}_{i=2}(p_i-1).$
  The tangent cone of this fast cycle  lies inside $\R_{\ge0}\hx_1\times \C^{r-1}_{x_2\dots x_r}$  and is of $\dim_\R \le r.$ The vanishing rate is $\frac{p_1}{p_{r+1}}.$

  Finally $V(f_p+f_{>p})$ is IMC \iff  $p_1=\cdots=p_n.$

 \item  (Semi-Brieskorn-Pham singularities,  the complete intersection case, $c\ge1$.) Take $f_{p_j}(\ux)=\sum_i a_{ji} x^{\quots{p_j}{\om_i}}_i.$
 Here $\quots{p_j}{\om_i}\in \N$ for all $j,i.$ Suppose the matrix $\{a_{ji}\}\in \Mat_{r\times N}(\C)$ is generic, to ensure that  the germ  $V(\uf_\up)$ is an ICIS. Take a higher order perturbation, $\uf_\up+\uf_{>\up}.$
  The assumptions of the corollary are implied by the genericity of the coefficients $\{a_{ji}\}.$   We get:    if $V(\uf_\up+\uf_{>\up})$ is IMC then $\om_1=\cdots=\om_n.$

\eee
\eex

\subsubsection{Perturbations of weighted-homogeneous hypersurfaces with non-isolated singularities} Take a hypersurface germ  $X=V(f_p+f_{>p})\sset (\C^{n+1},o),$ with $n\ge2$ and $f_p\in (x)^2.$ Here the germ $X_o=V(f_p)$ can have non-isolated singularity,  but we assume:
\bei
\item  $\dim_\C [X\cap \Sing (V(f_p))\smin V(x_1)]<\frac{n+1}{2}.$

\item Both $X $ and $X\cap V(x_1)$ (of $\dim_\C=n,n-1$) have isolated hypersurface singularities. Denote $\mu:=\mu(X\cap V(x_1)).$

\eei
\bcor\label{Thm.fast.cycles.on.perturb.weighted.homogen.hypersurf.sing}
\bee
\item If $\om_1\!<\!\om_n$ then $X$ has a fast cycle of homotopy type $\vee^\mu S^{n-1}$ and vanishing rate$\ge\!\frac{\om_{r+1}}{\om_1}.$

  \item Therefore, if $X$ is IMC then $\om_1=\cdots=\om_n.$
\eee
\ecor

\subsubsection{Infinity of exotic Lipschitz structures on $(\R^{2n},o)$}\label{Ex.Exotic.Spheres}
Take an isolated weighted-homogeneous hypersurface singularity $X\sset (\C^{n+1},o),$ with $n\ge3.$   In certain cases the links are exotic spheres, $\Link[X]\stackrel{homeo}{\approx}S^{2n-1},$ see e.g. \cite[pg.97-98]{Dimca}.
 E.g. $\Link[V(f)]$ is a topological sphere for $f(x)=\sum x^{p_i} $ when at least two of $p_i$'s are coprime to all the others.

Then one has the germ of a topological manifold, $X\stackrel{homeo}{\approx}(\R^{2n},o).$  However, for $p_1=\cdots=p_r<\cdots\le p_{n+1},$ with $r<n$, Corollary \ref{Thm.IMC.criterion.semi.weighted.homogen.germs} ensures a fast cycle on $X.$
  Hence, in these cases the homeomorphism $X\approx(\R^{2n},o)$ cannot be chosen inner-Lipschitz.
 Moreover we get infinity of distinct inner-Lipschitz cases.
 \bcor\label{Thm.Exotic.Structures}
 For each triple $(n,p,r),$ with $n\ge3,$ $p\ge2,$   $r<n,$ there exists a countable number of weighted-homogeneous hypersurface-germs $X\sset (\C^{n+1},o)$, each having $mult(X)=p,$ and weights $\om_1=\cdots=\om_r<\om_{r+1}\le\cdots\le \om_{n+1},$
 all homeomorphic to $(\R^{2n},o),$ but   of (pairwise) distinct inner-Lipschitz-types
 \ecor
This countable collection is reached in two ways:
\bei
\item either  via fast cycles of growing vanishing rates (by ensuring arbitrarily large $\frac{\om_{r+1}}{\om_1}$);
\item or via fast cycles of growing topological complexity, by ensuring arbitrarily large $\mu(X\cap V(x_1)).$
\eei

 The previously known examples of countably many exotic Lipschitz structures were constructed using separating sets, see Appendix to \cite{Birbrair.Fernandes.Neumann.Grandjean.O'Shea14}. Our example seems to be the first construction via fast cycles.

\subsubsection{Fast cycles on Newton-non-degenerate hypersurface germs}  Take a Newton-non-degenerate germ $X=V(f)\sset (\C^{n+1},o).$
 Assume  $f$ is convenient, i.e. its Newton diagram, $\Ga_f\sset \R^{n+1}_{\ge0},$ intersects all the coordinate axes.
 Below we consider only the top-dimensional faces,  $\si\sset \Ga_f.$
 Assume  $\Ga_f$ has at least two such faces.
  Each   $\si\sset \Ga_f$ has   weights, $\om_1,\dots,\om_{n+1}.$ These are the coordinates of the normal vector to $\si.$

\bcor\label{Thm.IMC.criterion.Newton.non.degen.germs}
 Suppose a face $\si\sset \Ga_f$ satisfies (after permuting the variables):
 \bee[i.]
 \item   $\om_1\le \cdots\le \om_{n+1},$  and $\si$ intersects the hyperplane $\{x_1=0\}\sset \R^{n+1}.$
 \item $\dim_\C [X\cap \Sing V(f_{|\si})\smin V(x_1)]<\frac{n+1}{2}.$
 \eee
If   $\om_1<\om_n$  then $X$ has a fast cycle of homotopy type $\vee^\mu S^{n-1},$ where $\mu=\mu(X\cap V(x_1)).$
\\
Therefore, if $X$ is IMC  then the lower $n$ weights of that face $\si$ coincide.
\ecor
\noindent Condition ii. is non-empty  as  $V(f_{|\si})$ can have non-isolated singularities, can be  non-reduced.
\bpr
  Restrict the function to this face, $f_{|\si},$ one gets the weighted-homogeneous hypersurface $V(f_{|\si})\sset  \C^{n+1} .$
   Then we present
 $f=f_{|\si}+f_{>\si}.$
 Thus the germ $X=V(f_{|\si}+f_{>\si})$ is a perturbation  of $X_o:=V(f_{|\si}).$
 The statement follows by  theorem \ref{Thm.IMC.criterion.via.weights}.
  Note that both $X$ and $X\cap V(x_1)$ have isolated singularities.
\epr

\bex\label{Ex.IMC.that.are.NND}
\bee[\bf i.] %\hspace{-0.5cm}
\item (Newton-non-degenerate surface germs, $n=2$) The condition $\dim_\C [V(f)\cap \Sing V(f_\si)]<\frac{n+1}{2}$ is satisfied (trivially) for each face.
 Thus the IMC assumption on $X$ implies:  for any face $\si\sset \Ga$  that intersects the coordinate plane of lowest weight,  the two lowest weights coincide.

\item  (Newton-non-degenerate three-fold germs, $n\=3$) The condition
\beq
\dim_\C[  V(f)\cap \Sing V(f_{|\si})\smin V(x_1)]\<\frac{n+1}{2}
\eeq
  means
  $\dim_\C [\Sing V(f_{|\si})\smin V(x_1)] <3.$ And the later means:   $V(f_{|\si})\smin V(x_1)$ is reduced. Which means: $f_{|\si}$ is not divisible by any of $x^2_2, x^2_3, x^2_4.$
 Equivalently: the distance of $\si$ from each of the coordinate hyperplanes $V(x_2),V(x_3),V(x_4)$ is at most one.
  If this is satisfied for $\si,$ and moreover $\si\!\cap \!\{x_1=0\}\!\neq\!\empty,$ and $V(f)\sset (\C^4,o)$ is IMC, then the three lowest weights of $\si$   coincide.
\eee
\eex

 \end{document}